# How Peircean was the "'Fregean' Revolution" in Logic?

Irving H. Anellis

Peirce Edition, Institute for American Thought
Indiana University – Purdue University at Indianapolis
Indianapolis, IN, USA
ianellis@iupui.edu

*Abstract*. The historiography of logic conceives of a Fregean revolution in which modern mathematical logic (also called symbolic logic) has replaced Aristotelian logic. The preeminent expositors of this conception are Jean van Heijenoort (1912–1986) and Donald Angus Gillies. The innovations and characteristics that comprise mathematical logic and distinguish it from Aristotelian logic, according to this conception, created *ex nihlo* by Gottlob Frege (1848–1925) in his *Begriffsschrift* of 1879, and with Bertrand Russell (1872–1970) as its chief This position likewise understands the algebraic logic of Augustus De Morgan (1806–1871), George Boole (1815–1864), Charles Sanders Peirce (1838–1914), and Ernst Schröder (1841–1902) as belonging to the Aristotelian tradition. The "Booleans" are understood, from this vantage point, to merely have rewritten Aristotelian syllogistic in algebraic guise.

The most detailed listing and elaboration of Frege's innovations, and the characteristics that distinguish mathematical logic from Aristotelian logic, were set forth by van Heijenoort. I consider each of the elements of van Heijenoort's list and note the extent to which Peirce had also developed each of these aspects of logic. I also consider the extent to which Peirce and Frege were aware of, and may have influenced, one another's logical writings.

§0. *The nature of the question in historical perspective.* Lest anyone be misled by the formulation of the question: "How Peircean was the "Fregean" Revolution in Logic?"; if we understand the question to inquire whether Peirce in some respect participated in the "Fregean revolution" or whether Peirce had in some wise influenced Frege or adherents of Frege's conception of logic, the unequivocal reply must be a decided: "No!" There is no evidence that Frege, at the time wrote his [*1879*] *Begriffsschrift* had even heard of Peirce, let alone read any of Peirce's writings in logic. More particularly, whatever Frege may have read by or about Peirce was by way of his subsequent interactions with Ernst Schröder that were opened by Schröder's [*1880*] review of the *Begriffsschrift*. What I have in mind in asking the question was whether there were elements in Peirce's logic or his conception of logic that have been identified as particularly characteristic of the "Fregean" conception of logic or novel contributions to logic which adherents of the historiographical conception of a "Fregean" revolution in logic have asserted were original to Frege, and which therefore distinguish the logic of Frege and the Fregeans as identifiably distinct from logic as it was previously known.

There are two ways of characterizing the essence of the "Fregean" revolution in logic. One, Jean van Heijenoort and Hans-Dieter Sluga among those adopting this view, asserts that Booleans are to be distinguished from Fregeans. This is a multi-faceted conception, the core of which is the notion that the Booleans saw logic as essentially algebraic, and regarded logic as a *calculus*, alongside of other algebras, whereas the Fregeans adopted a function-theoretic syntax and conceived of logic as preeminently a *language* which also happens to be a calculus. The other, led by Donald Angus Gillies, asserts that logic before Frege was Aristotelian. The criterion for the distinction between Aristotelians and Fregeans (or mathematical logicians) is whether the old subject-predicate syntax of proposition is adopted. Adherents of this line argue that the "Booleans" are also Aristotelians, their purpose being to simply rewrite Aristotelian propositions in symbolic form, to algebraicize Aristotle's syllogistic logic, to, in the words of William Stanley Jevons [*1864*, 3; *1890*, 5], clothe Aristotle in "mathematical dress."

AS editor of the very influential anthology *From Frege to Gödel: A Source Book in Mathematical Logic, 1879–1931* (hereafter FFTG) [van Heijenoort *1967*], historian of logic Jean van Heijenoort (1912–1986) did as much as anyone to canonize as historiographical truism the conception, initially propounded by Bertrand Russell (1872–1970), that modern logic began in 1879 with the publication of the *Begriffsschrift* [Frege *1879*] of Gottlob Frege (1848–1925). Van Heijenoort did this by relegating, as a minor sidelight in the history of logic, perhaps "interesting in itself" but of little historical impact, the tradition of algebraic logic of George Boole (1815–1864), Augustus De Morgan (1806–1871), Charles Sanders Peirce (1839–1914), William Stanley Jevons (1835–1882), and Ernst Schröder (1841–1902).

The first appearance of the *Begriffsschrift* prompted reviews in which the reviewers argued on the one hand that Frege's notational system was unwieldy (see, *e.g*. [Schröder *1880*, 87–90]), and, on the other, more critically, that it offered little or nothing new, and betrayed either an ignorance or disregard for the work of logicians from Boole forward. Schröder [*1880*, 83], for example, wrote that "In ersten Linie finde ich an der Schrift auszusetzen, dass dieselbe sich zu isolirt hinstellt und an Leistungen, welche in sachlich ganz verwandten Richtungen—namentlich von Boole gemacht sind, nicht nur keinen ernstlichen Anschluss sucht, sondern dieselben gänzlich unberücksichtigt lässt." Frege's *Begriffsschrift* is not nearly so essentially different from Boole's formal language as is claimed for it, Schröder [*1880*, 83] adds, declaring [Schröder *1880*, 84] that one could even call the *Begriffsschrift* an "Umschreibung", a paraphrase, of Boole's formal language.

In many respects, the attitudes of Frege and Edmund Husserl (1859–1938) toward algebraic logic were even more strongly negative. We recall, for example, the chastisement by Schröder's student Andreas Heinrich Voigt (1860–1940) of Husserl's assertion in "Der Folgerungscalcul und die Inhaltslogik" [Husserl *1891a*, 171]—"nicht…in dem gewöhnlichen Sinne der Logik"—and the review [Husserl *1891b*, 246–247] of the first volume [Schröder *1890*] of Schröder's *Vorlesungen über die Algebra der Logik* that algebraic logic is not logic [Voigt *1892*], and Frege's ire at Husserl [*1891b*, 243] for regarding Schröder, rather than Frege, as the first in Germany to attempt in symbolic logic, and indeed the first in Germany to attempt to develop a "full-scale" extensional logic.[1] Not only that; Voigt in "Zum Calcul der Inhaltslogik. Erwiderung auf Herrn Husserls Artikel" [*1893*, 506] pointed out that much of what Husserl claimed as original for his logic was already to be found in Frege and Peirce, pointing in particular to Peirce's [*1880*] "On the Algebra of Logic".[2] Responding to remarks by Schröder to his doctoral dissertation, Voigt [*1893*, 506–507] informs readers that in revising his dissertation, he wrote, in part, dealing with "die Logik der Gattungen (des Inhalts)":

> „Gewöhnlich sind die Bearbeiter der algebraischen Logik von der Anschauung ausgegangen, dass alle Begriffe als Summen von Individuen, d. h. als Classen anzusehen seien, und man hält daher in Folge dessen häufig diese Anschauungsweise für eine der algebraischen Logik wesentlich, über die sie auch nicht hinauskönne. Dass dieses keineswegs der Fall, dass sie sogut wie die alte Logik auch eine Logik des Inhalts sein kann, hat, soviel ich weiss, zuerst Herr FREGE (Begriffsschrift, Halle a. S. 1879), dann besonders Herr PEIRCE (a. a. O. 1880) gezeigt, und wenn auch in der Begründung einiger Principien bei Peirce noch eine kleine Lücke ist, so hat es doch keine Bedenken, diese Principien axiomatisch gelten zu lassen, u. s. w."

He then writes [Voigt *1893*, 507]:

> In diesem Stück meiner Dissertation sind schon die zwei Schriftsteller erwähnt, welche einen Logikcalul unabhängig von Classenbeziehungen begründet haben, und von deren Herr HUSSERL weinigstens FREGE hätte kennen können, wenn ihn auch die Hauptarbeit PEIRCE'S im American Journal of Mathematics, Vol. III, nicht zugänglich gewesen wäre. FREGE hat einen leider in der Form sehr unbeholfenen, im Wesen aber mit dem SCHRÖDER'schen und jedem anderen Calcul übereinstimmenden Calcul geschaffen. Ueberhaupt steht es wohl von vorherein fest, dass jeder logische Calcul, wie er auch begründet werden mag, nothwending mit den bestehenden Calculen in Wesentlichen übereinstimmen muss.

Both Voigt and Husserl argue that their own respective logical systems are both a contentual logic (Inhaltslogik) and a deductive system, hence both extensional and intensional, and hence, as deductive, a calculus and, as contentual, a language. Husserl understands Schröder's algebra of logic, however, to be merely a calculus, concerned, he asserts [Husserl *1891b*, 244] exclusively with deduction, denying that Schröder's manifolds or sets (the *Mannigfaltigkeiten*, *i.e*. Schröder's classes) are legitimately extensions. Husserl [*1891b*, 246] rhetorically asks: "Ist aber Rechnen ein Schließen?" His answer [Husserl *1891b*, 246] immediately follows: "Keineswegs. Das Rechnen ist ein blindes Verfahren mit Symbolen nach mechanisch-reproducierten Regeln der Umwandlung und Umsetzung von Zeichen des jeweilgen Algorithmus." And this is precisely what we find in Schröder's *Algebra der Logik*, and nothing else. Husserl [*1891b*, 258] expands upon his assertion, explaining why Schröder is confused and incorrect in thinking his algebra is a logic which is a language rather than a mere calculus, *viz*.:

> Es ist nicht richtig, daβ die ›exacte‹ Logik nichts anderes ist, als ein Logik auf Grund einen neuen Sprache. Sie ist...überhaupt keine Logik, sondern ein speciellen logischen Zwecken dienender Calcul, und so ist denn die Rede von einer »Darstellung der Logik als einer Algebra« eine ganz unpassende,

Schröder's error, in Husserl's estimation, was to confuse or conflate a language with an algorithm, and hence fail to differentiate between a language and a calculus. He defines a calculus [Husserl *1891b*, 265] as nothing other than a system of formulas, entirely in the manner of externally based conclusions—"nichts Anderes als ein System des formellen, d.i. rein auf die Art der Aeuβer sich gründenden Schlieβens."

At the same time, Husserl in his "Antwort auf die vorstehende, Erwiderung' der Herrn Voigt" denies that Frege's *Begriffsschrift* is a calculus in the exact meaning of the word [Husserl *1893*, 510], and neither does Peirce have more than a calculus, although he credits Husserl with at least having the concept of a content logic [Husserl *1893*, 510]. Nevertheless, not until afterwards, in his anti-psychologistic *Logische Untersuchungen* [Husserl *1900-1901*] which served as the founding document of his phenomenology, did Husserl see his study of logic as the establishment of formal logic as a *characteristica universalis*.

In response to Husserl's [*1891a*, 176] assertion:

> Vertieft man sich in die verschiedenen Versuche, die Kunst der reinen Folgerungen auf eine calculierende Tevhnik zu bringen, so merkt man wesentlich Unterschiede gegenüber der Verfahrungsweisen der alten Logik. ...Und mit vollen Rechte, wofern sie nur den Anspruch nicht mehr erhebt, statt einer blossen Technik des Folgerns, eine Logik derselben zu bedeuten, in which Boole, Jevons, Peirce and Schröder are identified by Husserl [*1891a*, 177] as developers of an algorithmic calculus of inference rather than a true logic, Voigt [*1892*, 292] asserts that the algebra of logic is just as fully a logic as the older—Aristotelian—logic, having the same content and goals, but more exact and reliable,[3] he takes aim at Husserl's claim [Husserl *1891a*, 176] that algebraic logic is not a logic, but a calculus, or, in Husserl's words, only a symbolic technique; "dass die Algebra der Logik keine Logik, sondern nur ein Calcul der Logik, eine mechanische Methode nicht der logischen Denkens, sondern sich logisches Denken zu ersparen, sei" [Voigt *1892,* 295]. Voigt notes, Husserl confuses deductive inference with mental operations. Husserl denies that the algebra of logic is deductive, arguing that it cannot examine its own inference rules, since it is limited to concepts. Voigt [*1892*, 310] replies by remarking that, if Husserl is correct, then neither is syllogistic logic deductive, and he then defines deductive logic as concerned with the relations between concepts and judgments and notes that the second volume of Schröder's *Vorlesungen…* [Schröder *1891*] indeed introduces judgments. He demonstrates how to write equations in Schröder's system that are equivalent to categorical syllogisms, and presents [Voigt *1892*, 313*ff*.] in Schröder's notation the Aristotle's logical Principles of Identity, Non-contradiction, and Excluded Middle, along with the laws of distribution and other algebraic laws to demonstrate that the algebra of logic, composed of both a logic of judgments and a logic of concepts indeed is a deductive logic.

This claim that Schröder's algebra of logic is not a logic also found its echo in Frege's review of the first volume of Schröder's [*1890*] *Vorlesungen über die Algebra der Logik* when he wrote [Frege *1895*, 452] that: "Alles dies ist sehr anschaulich, unbezweifelbar; nur schade: es ist unfruchtbar, und es ist keine Logik." The reason, again, is that Schröder's algebra does not deal with relations between classes. He goes so far as to deny even that it is a deductive logic or a logic of inferences. He says of Schröder's algebra of logic [Frege *1895*, 453] that it is merely a calculus, in particular a *Gebietkalkul*, a domain calculus, restricted to a Boolean universe of discourse; and only when it is possible to express thoughts in general by dealing with relations between classes does one attain a logic—"nur dadurch [allgemein Gedanken auszudrücken, indem man Beziehungen zwischen Klassen angiebt]; nur dadurch gelangt man zu einer Logik."

In asserting that the algebraic logicians present logic as a calculus, but not logic as a language, van Heijenoort is, in effect arguing the position taken by Frege and Husserl with respect to Schröder's algebra of logic, that it is a mere calculus, not truly or fully a logic. It is the establishment of logic as a language that, for Frege and for van Heijenoort, constitute the essential difference between the Booleans or algebraic logicians and the quantification-theoretical mathematical logicians, and encapsulates and establishes the essence of the Fregean revolution in logic.

Russell was one of the most enthusiastic early supporters of Frege and contributed significantly to the conception of Frege as the originator of modern mathematical logic, although he never explicitly employed the specific term "Fregean revolution". In his recollections, he states that many of the ideas that he thought he himself originated, he later discovered had already been first formulated by Frege, and some others were due to Giuseppe Peano (1858–1932) or the inspiration of Peano.

The conception of a Fregean revolution was further disseminated and enhanced in the mid-1920s thanks to Paul Ferdinand Linke (1876–1955), Frege's friend and colleague at Jena helped formulate the concept of a "Fregean revolution" in logic, when he wrote [Linke *1926*, 226–227] that:

> …the great reformation in logic…originated in Germany at the beginning of the present century…was very closely connected, at least at the outset, with mathematical logic. For at bottom it was but a continuation of ideas first expressed by the Jena mathematician, Gottlob Frege. This prominent investigator has been acclaimed by Bertrand Russell to be the first thinker who correctly understood the nature of numbers. And thus Frege played an important role in…mathematical logic, among whose founders he must be counted.

Russell's extant notes and unpublished writing demonstrate that significant parts of logic that he claimed to have been the first to discover were already present in the logical writings of Charles Peirce and Ernst Schröder. With regard to Russell's claim, to having invented the logic of relations, he was later obliged to reluctantly admit that Peirce and Schröder had already "treated" of the subject, so that, in light of his own work, it was unnecessary to "go through" them.

We also find that Bertrand Russell (1872–1970) not only had read Peirce's "On the Algebra of Logic" [Peirce *1880*] and "On the Algebra of Logic: A Contribution to the Philosophy of Notation" [Peirce *1885*] and the first volume of Schröder's *Vorlesungen über die Algebra der Logik* [Schröder *1890*] earlier than his statements suggest,[4] and had known the work and many results even earlier, in the writing of his teacher Alfred North Whitehead (1861–1947), as early as 1898, if not earlier, indeed when reading the galley proofs of Whitehead's *Treatise of Universal Algebra* [Whitehead *1898*], the whole of Book II, out of seven of which was devoted to the "Algebra of Symbolic Logic", came across references again in Peano, and was being warned by Louis Couturat (1868–1914) not to short-change the work of the algebraic logicians. (For specific examples and details, including references and related issues, see [Anellis *1990/1991*; *1995*; *2004-2005*; *2011*], [Hawkins *1997*].) There is of course also published evidence of Russell at the very least being aware that Peirce and Schröder worked in the logic of relatives, by the occasional mention, however denigratory and haphazard, in his *Principles of Mathematics*.[5]

For the greater part, my approach is to reorganize what is—or should be—already known about Peirce's contributions to logic, in order to determine whether, and if so, to what extent, Peirce's work falls within the parameters of van Heijenoort's conception of the Fregean revolution and definition of mathematical logic, as particularized by the seven properties or conditions which van Heijenoort listed as characterizing the "Fregean revolution" and defining "mathematical logic". I am far less concerned here with analyzing or evaluating van Heijenoort's characterization and the criterion which he lists as constituting Frege's revolution. The one exception in my rendition of Peirce's work is that I cite material to establish beyond any doubt that Peirce had developed truth table matrices well in advance of the earliest examples of these, identified by John Shosky [*1997*] as jointly attributable to Bertrand Russell and Ludwig Wittgenstein (1889–1951) and dating to 1912.

§1. *The defining characteristics of the "Fregean revolution"*. What historiography of logic calls the "Fregean revolution" was articulated in detail by Jean van Heijenoort.

In his anthology *From Frege to Gödel*, first published in 1967, and which historiography of logic has for long taken as embracing all of the significant work in mathematical logic, van Heijenoort [*1967*, vi] described Frege's *Begriffsschrift* of 1879 as of significance for the significance of logic, comparable, if at all, only with Aristotle's *Prior Analytics*, as opening "a great epoch in the history of logic…." Van Heijenoort listed those properties that he considered as characterizing Frege's achievements and that distinguishes modern mathematical logic from Aristotelian logic. These characteristics are such that the algebraic logic of Boole, De Morgan, Jevons, Peirce, Schröder and their adherents is regarded as falling outside the realm of modern mathematical logic, or, more precisely, are not properly considered as included within the purview of modern mathematical logic as formulated and developed by, or within, the "Fregean revolution".

In his posthumously published "Historical Development of Modern Logic", originally composed in 1974, he makes the point more forcefully still of the singular and unmatched significance of Frege and his *Begriffsschrift* book-let of a mere 88 pages; he began this essay with the unequivocal and unconditional declaration [van Heijenoort *1992*, 242] that: "Modern logic began in 1879, the year in which Gottlob Frege (1848–1925) published his *Begriffsschrift*. He then goes on, to explain [van Heijenoort *1992*, 242] that:

> In less than ninety pages this booklet presented a number of discoveries that changed the face of logic. The central achievement of the work is the theory of quantification; but this could not be obtained till the traditional decomposition of the proposition into subject and predicate had been replaced by its analysis into function and argument(s). A preliminary accomplishment was the propositional calculus, with a truth-functional definition of the connectives, including the conditional. Of cardinal importance was the realization that, if circularity is to be avoided, logical derivations are to be *formal,* that is, have to proceed according to rules that are devoid of any intuitive logical force but simply refer to the typographical form of the expressions; thus the notion of formal system made its appearance. The rules of quantification theory, as we know them today, were then introduced. The last part of the book belongs to the *foundations of mathematics,* rather than to logic, and presents a logical definition of the notion of mathematical sequence. Frege's contribution marks one of the sharpest breaks that ever occurred in the development of a science.

We cannot help but notice the significant gap in the choices of material included in FFTG—all of the algebraic logicians are absent, not only the work by De Morgan and Boole, some of which admittedly appeared in the late 1840s and early 1850s, for example Boole's [*1847*] *The Mathematical Analysis of Logic* and [*1854*] *An Investigation of the Laws of Thought,* De Morgan's [*1847*] *Formal Logic,* originating algebraic logic and the logic of relations, and Jevons's "de-mathematicizing" modifications of Boole's logical system in his [*1864*] *Pure Logic or the Logic of Quality apart from Quantity* and [*1869*] *The Substitution of Similars*, not only the first and second editions of John Venn's *Symbolic Logic* [Venn *1881*; *1894*] which, along with Jevons's logic textbooks, chiefly his [*1874*] *The Principles of Science, a Treatise on Logic and Scientific Method* which went into its fifth edition in 1887, were particularly influential in the period from 1880 through 1920 in disseminating algebraic logic and the logic of relations, but even for work by Peirce and Schröder that also appeared in the years which this anthology, an anthology purporting to completeness, includes, and even despite the fact that Frege and his work is virtually unmentioned in any of the other selections, whereas many of the work included do refer back, often explicitly, to contributions in logic by Peirce and Schröder. Boole's and De Morgan's work in particular served as the starting point for the work of Peirce and Schröder.

The work of these algebraic logicians is excluded because, in van Heijenoort's estimation, and in that of the majority historians and philosophers—almost all of whom have since at least the 1920s, accepted this judgment, the work of the algebraic logicians falls outside of the Fregean tradition. It was, however, far from universally acknowledged during the crucial period between 1880 through the early 1920s, that either Whitehead and Russell's *Principia Mathematica* nor any of the major efforts by Frege, was the unchallenged standard for what mathematical logic was or ought to look like.[6]

Van Heijenoort made the distinction one primarily between algebraic logicians, most notably Boole, De Morgan, Peirce, and Schröder, and logicians who worked in quantification theory, first of all Frege, and with Russell as his most notable follower. For that, the logic that Frege created, as distinct from *algebraic* logic, was regarded as *mathematical* logic. ([Anellis & Houser *1991*] explore the historical background for the neglect which algebraic logic endured with the rise of the "modern mathematical" logic of Frege and Russell.)

Hans-Dieter Sluga, following van Heijenoort's distinction between followers of Boole and followers of Frege, labels the algebraic logicians "Booleans" after George Boole, thus distinguishes between the "Fregeans", the most important member of this group being Bertrand Russell, and the "Booleans", which includes not only, of course, Boole and his contemporary Augustus De Morgan, but logicians such as Peirce and Schröder who combined, refined, and further developed the algebraic logic and logic of relations established by Boole and De Morgan (see [Sluga *1987*]).

In the last two decades of the nineteenth century and first two decades of the twentieth century, it was, however, still problematic whether the new Frege-Russell conception of mathematical logic or the classical Boole-Schröder calculus would come to the fore. It was also open to question during that period whether Russell (and by implication Frege) offered anything new and different from what the algebraic logicians offered, or whether, indeed, Russell's work was not just a continuation of the work of Boole, Peirce, and Schröder Peano (see [Anellis *2004-2005*; *2011*]), for example, such of regarded Russell's work as "On Cardinal Numbers" [Whitehead *1902*], §III of which was actually written solely by Russell) and "Sur la logique des relations des applications à la théorie des séries" [Russell *1901a*] as "filling a gap" between his own work and that of Boole, Peirce, and Schröder (see [Kennedy *1975*, 206]).[7] Through the *fin de siècle*, logicians for the most part understood Russell to be transcribing into Peanesque notation Cantorian set theory and the logic of relations of Peirce and Schröder.

Bertrand Russell, in addition to the strong and well-known influence which Giuseppe Peano had on him, was a staunch advocate, and indeed one of the earliest promoters, of the conception of a "Fregean revolution" in logic, although he himself never explicitly employed the term itself. Nevertheless, we have such pronouncements, for example in his manuscript on "Recent Italian Work on the Foundations of Mathematics" of 1901 in which he contrasts the conception of the algebraic logicians with that of Hugh MacColl (1837–1909) and Gottlob Frege, by writing that (see [Russell *1993*, 353]):

> Formal Logic is concerned in the main and primarily with the relation of implication between propositions. What this relation is, it is impossible to define: in all accounts of Peano's logic it appears as one of his indefinables. It has been one of the bad effects of the analogy with ordinary Algebra that most formal logicians (with the exception of Frege and Mr. MacColl) have shown more interest in logical equations than in implication.

§2. *The characteristics of modern mathematical logic as defined and delimited by the "Fregean revolution"*. In elaborating the distinguishing characteristics of mathematical logic and, equivalently, enumerating the innovations which Frege—allegedly—wrought to create mathematical logic, van Heijenoort (in "Logic as Calculus and Logic as Language" [van Heijenoort *1976b*, 324]) listed:

1. a propositional calculus with a truth-functional definition of connectives, especially the conditional;
2. decomposition of propositions into function and argument instead of into subject and predicate;
3. a quantification theory, based on a system of axioms and inference rules; and
4. definitions of *infinite sequence* and *natural number* in terms of logical notions (*i.e*. the logicization of mathematics).

In addition, Frege, according to van Heijenoort and adherents of the historiographical conception of a "Fregean revolution":

5. presented and clarified the concept of *formal system*; and
6. made possible, and gave, a use of logic for philosophical investigations (especially for philosophy of language).

Moreover, in the undated, unpublished manuscript notes "On the Frege-Russell Definition of Number",[8] van Heijenoort claimed that Bertrand Russell was the first to introduce a means for

7. separating singular propositions, such as "Socrates is mortal" from universal propositions such as "All Greeks are mortal."

In the "Historical Note" to the fourth edition of his *Methods of Logic* [Quine *1982*, 89] Willard Van Orman Quine (1908–2000) asserts that "Frege, in 1879, was the first to axiomatize the logic of truth functions and to state formal rules of inference." Similar remarks are scattered throughout his textbook. He does, however, give Pierce credit [Quine *1982*, 39]—albeit along with Frege and Schröder—for the "pattern of reasoning that the truth table tabulates."

Defenders of the concept of a Fregean revolution" count Peirce and Schröder among the "Booleans" rather than among the "Fregeans". Yet, judging the "Fregean revolution" by the (seven) supposedly defining characteristics of modern mathematical logic, we should include Peirce as one of its foremost participants, if not one of its initiators and leaders. At the very least, we should count Peirce and Schröder among the "Fregeans" rather than the "Booleans" where they are ordinarily relegated and typically have been dismissed by such historians as van Heijenoort as largely, if not entirely, irrelevant to the history of modern mathematical logic, which is "Fregean".

Donald Gillies is perhaps the leading contemporary adherent of the conception of the "Fregean revolution", and he has emphasized in particular the nature of the revolution as a replacement of the ancient Aristotelian paradigm of logic by the Fregean paradigm. The centerpiece of this shift is the replacement of the subject-predicate syntax of Aristotelian propositions by the function-argument syntax offered by Frege (*i.e*. van Heijenoort's second criterion). They adhere to the subject-predicate structure for propositions.

Whereas van Heijenoort and Quine stressed in particular the third of the defining characteristics of Fregean or modern mathematical logic, the development of a quantification theory, Gillies [*1992*] argues in particular that Boole and the algebraic logicians belong to the Aristotelian paradigm, since, he explains, they understood themselves to be developing that part of Leibniz's project for establishing a *mathesis universalis* by devising an arithmeticization or algebraicization of Aristotle's categorical propositions and therefore of Aristotelian syllogistic logic, and therefore retain, despite the innovations in symbolic notation that they devised, the subject-predicate analysis of propositions.

What follows is a quick survey of Peirce's work in logic, devoting attention to Peirce's contributions to all seven of the characteristics that distinguish the Fregean from the Aristotelian or Boolean paradigms. While concentrating somewhat on the third, which most defenders of the conception of a "Fregean revolution" count as the single most crucial of those defining characteristics. The replacement of the subject-predicate syntax with the function-argument syntax is ordinarily accounted of supreme importance, in particular by those who argue that the algebraic logic of the "Booleans" is just the symbolization, in algebraic guise, of Aristotelian logic. But the question of the nature of the quantification theory of Peirce, Mitchell, and Schröder as compared with that of Frege and Russell is tied up with the ways in which quantification is handled.

The details of the comparison and the mutual translatability of the two systems is better left for another discussion. Suffice it here to say that Norbert Wiener (1894–1964), who was deeply influenced by Royce and was the student of Edward Vermilye Huntington (1874–1952), a mathematician and logician who had corresponded with Peirce, dealt with the technicalities in detail in his doctoral thesis for Harvard University of 1913, *A Comparison Between the Treatment of the Algebra of Relatives by Schroeder and that by Whitehead and Russell* [Wiener *1913*], and concluded that there is nothing that can be said in the *Principia Mathematica* (1910-13) of Whitehead and Russell that cannot, with equal facility, be said in the Peirce-Schröder calculus, as presented in Schröder's *Vorlesungen über die Algebra der Logik* [Schröder *1890-1905*]. ([Grattan-Guinness *1975*] is a discussion of Wiener's thesis.) Indeed, Brady [*2000*, 12] essentially asserts that Wiener accused Russell of plagiarizing Schröder, asserting, without giving specific references, that Wiener [*1913*] presents "convincing evidence to show" that Russell "lifted his treatment of binary relations in *Principia Mathematica* almost entirely from Schröder's *Algebra der Logik*, with a simple change of notation and without attribution." In his doctoral thesis, Wiener had remarked that "Peirce developed an algebra of relatives, which Schröder extended…." Russell could hardly have missed that assertion; but it was in direct contradiction to one of Russell's own self-professed claims to have devised the calculus of relations on his own. Russell complained in reply that Wiener considered only "the more conventional parts of *Principia Mathematica*" (see [Grattan-Guinness *1975*, 130]). Thereafter, Wiener received a traveling fellowship from Harvard that took him to Cambridge from June 1913 to April 1914 and saw him enrolled in two of Russell's courses, one of which was a reading course on *Principia Mathematica*, and in a mathematics course taught by G. H. Hardy. They met repeatedly between 26 August and 9 September 1913 to discuss a program of study for Wiener. In these discussions, Russell endeavored to convince Wiener of the greater perspicacity of the *Principia* logic. Within this context they discussed Frege's conception of the *Werthverlauf* (course-of-values) and Russell's concept of propositional functions. Frege's [*1893*] *Grundgesetze der Arithmetik*, especially [Frege *1893*, §11], where Frege's function $\backslash\xi$ replaces the definite article, such that, for example, $\backslash\xi$(positive $\sqrt{2}$) represents the concept which is the proper name of the positive square root of 2 when the value of the function $\backslash\xi$ is the positive square root of 2, and to Peano's [*1897*] "Studii di logica matematica", in which Peano first considered the role of "the", the possibility its elimination from his logical system; whether it can be eliminated from mathematical logic, and if so, how. In the course of these discussions, Russell raised this issue with Norbert Wiener (see [Grattan-Guinness *1975*, 110]), explaining that:

> There is need of a notation for "the". What is alleged does not enable you to put "0 = etc. Df". It was a discussion on this very point between Schröder and Peano in 1900 at Paris that first led me to think Peano superior.

After examining and comparing the logic of *Principia* with the logic of Schröder's *Vorlesungen über die Algebra der Logic*, Wiener developed his simplification of the logic of relations, in a very brief paper titled "A Simplification of the Logic of Relations" [Wiener *1914*] in which the theory of relations was reduced to the theory of classes by the device of presenting a definition, borrowed from Schröder, of ordered pairs (which, in Russell's notation, reads $\langle x, y\rangle = \{\{\{x\}, \varnothing\}, \{y\}\}\}$), in which a relation is a class of ordered couples. It was sufficient to prove that $\langle a, b\rangle = \langle c, d\rangle$ implies that $a = b$ and $c = d$ for this definition to hold.

In a consideration that would later find its echo in Wiener's comparison, Voigt argued that Frege's and Schröder's systems are equivalent.

With that in mind, I want to focus attention on the question of quantification theory, without ignoring the other technical points.

1. *Peirce's propositional calculus with a truth-functional definition of connectives, especially the conditional*:

Peirce's contributions to propositional logic have been studied widely. Attention has ranged from Anthony Norman Prior's (1914–1969) expression and systematization of Peirce's axioms for the propositional calculus [Prior *1958*] to Randall R. Dipert's survey of Peirce's propositional logic [Dipert *1981*]. It should be clear that Peirce indeed developed a propositional calculus, which he named the "icon of the first kind" in his publication "On the Algebra of Logic: A Contribution to the Philosophy of Notation" [Peirce *1885*].

In an undated, untitled, two-page manuscript designated "Dyadic Value System" (listed in the Robin catalog as MS #6; [Peirce *n.d.(a)*][9]), Peirce asserts that the simplest of value systems serves as the foundation for mathematics and, indeed, for all reasoning, because the purpose of reasoning is to establish the truth or falsity of our beliefs, and the relationship between truth and falsity is precisely that of a dyadic value system, writing specifically that "the the whole relationship between the values," 0 and 1 in what he calls a cyclical system "may be summed up in two propositions, first, that there are different values" and "second, that there is no third value." He goes on to says that: "With this simplest of all value-systems mathematics must begin. Nay, all reasoning must and does begin with it. For to reason is to consider whether ideas are <u>true</u> or <u>false</u>." At the end of the first page and the beginning of the second, he mentions the principles of Contradiction and Excluded Middle as central. In a fragmented manuscript on "Reason's Rules" of *circa* 1902 [Peirce *ca. 1902*], he examines how truth and falsehood relate to propositions.

Consider the formula $[(\sim c \supset a) \supset (\sim a \supset c)] \supset \{(\sim c \supset a) \supset [(c \supset a) \supset a]\}$ of the familiar propositional calculus of today. Substituting Peirce's hook or "claw" of illation ($-\!\!<$) or Schröder's subsumption (€) for the "horseshoe ($\supset$) and Peirce's over-bar or Schröder's negation prime for the tilde of negation suffices to yield the same formula in the classical Peirce-Schröder calculus; thus:

Peano-Russell: $[(\sim c \supset a) \supset (\sim a \supset c)] \supset \{(\sim c \supset a) \supset [(c \supset a) \supset a]\}$
Peirce: $[(\bar{c} -\!\!< a) -\!\!< (\bar{a} -\!\!< c)] -\!\!< \{(\bar{c} -\!\!< a) -\!\!< [(c -\!\!< a) -\!\!< a]\}$
Schröder: [(*c*′ € *a*) € (*a* ′ € *c*)] € {(*c*′ € *a*) € [(*c* € *a*) € *a*]}

One of Husserl's arguments against the claim that the algebra of logic generally, including the systems of Boole and Peirce, and Schröder's system in particular [Husserl *1891b*, 267–272] is the employment of 1 and 0 rather than truth-values *true* and *false*. Certainly neither Boole nor Peirce had not been averse to employing Boolean values (occasionally even using '∞' for universes of discourses of indefinite or infinite cardinality) in analyzing the truth of propositions. Husserl, however, made it a significant condition of his determination of logic as a calculus, as opposed to logic as a language, that truth-values be made manifest, and not represented by numerical values, and he tied this to the mental representation which languages serve.

In the manuscript "On the Algebraic Principles of Formal Logic" written in the autumn of 1879—the very year in which Frege's *Begriffsschrift* appeared, Peirce (see [Peirce *1989*, 23]) explicitly identified his "claw" as the "copula of inclusion" and defined material implication or logical inference, illation, as

1st, $A -\!\!< A$, whatever $A$ may be.
2nd If $A -\!\!< B$, and $B -\!\!< C$, then $A -\!\!< C$.

From there he immediately connected his definition with truth-functional logic, by asserting [Peirce *1989*, 23] that

> This definition is sufficient for the purposes of formal logic, although it does not distinguish between the relation of inclusion and its converse. Were it desirable thus to distinguish, it would be sufficient to add that the real truth or falsity of $A -\!\!< B$, supposes the existence of $A$.

The following year, Peirce continued along this route: in "The Algebra of Logic" of 1880 [Peirce *1880*, 21; *1989*, 170],

$$A -\!\!< B$$

is explicitly defined as "*A* implies *B*", and

$$A \mathrel{\overline{-\!\!<}} B$$

defines "*A* does not imply *B*." Moreover, we are able to distinguish universal and particular propositions, affirmative and negative, according to the following scheme:

| | | | |
|---|---|---|---|
| A. | $a -\!\!< b$ | All *A* are *B* | (universal affirmative) |
| E. | $a -\!\!< \bar{b}$ | No *A* is *B* | (universal negative) |
| I. | $\breve{a} -\!\!< b$ | Some *A* is *B* | (particular affirmative) |
| O. | $\breve{a} -\!\!< \bar{b}$ | Some *A* is not *B* | (particular negative) |

In 1883 and 1884, in preparing preliminary versions for his article "On the Algebra of Logic: A Contribution to the Philosophy of Notation" [Peirce *1885*], Peirce develops in increasing detail the truth functional analysis of the conditional and presents what we would today recognize as the indirect or abbreviated truth table.

In the undated manuscript "Chapter III. Development of the Notation" [Peirce *n.d.(c)*], composed *circa* 1883-1884, Peirce undertook an explanation of material implication (without, however, explicitly terming it such), and making it evident that what he has in mind is what we would today recognize as propositional logic, asserting that letters represent assertions, and exploring the conditions in which inferences are valid or not, i.e., undertaking to "develope [*sic*] a calculus by which, from certain assertions assumed as premises, we are to deduce others, as conclusions." He explains, further, that we need to know, given the truth of one assertion, how to determine the truth of the other.

And in 1885, in "On the Algebra of Logic: A Contribution to the Philosophy of Notation" [Peirce *1885*], Peirce sought to redefine categoricals as hypotheticals and presented a propositional logic, which he called *icon of the first kind*. Here, Peirce [*1885*, 188–190], Peirce considered the *consequentia*, and introduces inference rules, in particular *modus ponens*, the "icon of the second kind" [Peirce *1885*, 188], transitivity of the copula or "icon of the third kind" [Peirce *1885*, 188–189], and *modus tollens*, or "icon of the fourth kind" [Peirce *1885*, 189].

In the manuscript fragment "Algebra of Logic (Second Paper)" written in the summer of 1884, Peirce (see [Peirce *1986*, 111–115]) reiterated his definition of 1880, and explained in greater detail there [Peirce *1986*, 112] that: "In order to say "If it is $a$ it is $b$," let us write $a -\!\!< b$. The formulae relating to the symbol $-\!\!<$ constitute what I have called the algebra of the copula…. The proposition $a -\!\!< b$ is to be understood as true if either $a$ is false or $b$ is true, and is only false if $a$ is true while $b$ is false."

It was at this stage that Peirce undertook the truth-functional analysis of propositions and of proofs, and also introduced specific truth-functional considerations, saying that, for **v** is the symbol for "true" (*verum*) and **f** the symbol for false (*falsum*), the propositions $\mathbf{f} -\!\!< a$ and $a -\!\!< \mathbf{v}$ are true, and either one or the other of $\mathbf{v} -\!\!< a$ or $a -\!\!< \mathbf{f}$ are true, depending upon the truth or falsity of $a$, and going on to further analyze the truth-functional properties of the "claw" or "hook".

In Peirce's conception, as found in his "Description of a Notation for the Logic of Relatives" of 1870, then Aristotelian syllogism becomes a hypothetical proposition, with material implication as its main connective; he writes [Peirce *1870*, 518] *Barbara* as:[10]

$$\begin{array}{r} \text{If } x -\!\!< y, \\ \text{and } y -\!\!< z, \\ \text{then } x -\!\!< z. \end{array}$$

In Frege's *Begriffsschrift* notation of 1879, §6, this same argument would be rendered as:

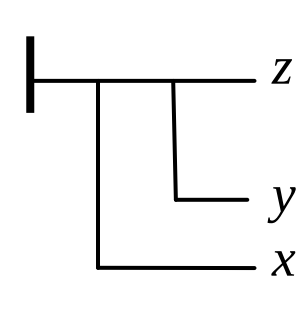

In the familiar Peano-Russell notation, this is just

$$(x \supset y) \cdot (y \supset z)] \supset (x \supset y).$$

Schröder, ironically, even complained about what he took to be Peirce's (and Hugh MacColl's) efforts to base logic on the propositional calculus, which he called the "MacColl-Peircean propositional logic."

Frege [*1895*, 434] recognized that implication was central to the logical systems of Peirce and Schröder (who employed '€', or Subsumption, in lieu of Peirce's '—<'), although criticizing them for employing the same symbol for class inclusion (or ordering) and implication, and thus for allegedly failing distinguish between these; class and set are in Schröder, he says [Frege *1895*, 435] "eingemischt", and which, in actuality, is just the part-whole relation. Thus he writes [Frege *1895*, 434]:

> Was Herr Schröder ‚Einordnung' oder ‚Subsumption' nennt, ist hier eigentlich nichts Anderes als die Beziehung des Teiles zum Ganzen mit der Erweiterung, dass jedes Ganze als seiner selbst betrachtet werden soll.

Frege [*1895*, 441–442] thus wants Schröder to distinguish between the s*ubter*-Beziehung, the class-inclusion relation, which is effectively implication, referencing [Frege *1895*, 442*n*.] in this regard Peano's [*1894*, §6] 'Ɔ', and the *sub*-Beziehung, or set membership relation, referencing [Frege *1895*, 442*n*.] Peano's [*1894*, §6] '$\in$'.

John Shosky [*1997*] distinguished between the truth-table *technique* or method on the one hand and the truth-table *device* on the other, the former amounting to a truth-functional analysis of propositions or arguments, the latter being in essence the presentation of truth-functional analysis in a tabular, or matrix, array. On this basis he argued that truth tables first appeared on paper in recognizable form around 1912, composed in the hand of either Ludwig Wittgenstein, with an addition by Bertrand Russell, on the verso of a typescript of a lecture by Russell on logical atomism, and thus earlier than its appearance in Wittgenstein's *Tractatus Logico-philosophicus* [Wittgenstein *1922*, Prop. 4.31] or the work of Emil Leon Post (1897–1954) and Jan Łukasiewicz (1878–1956) in 1920.[11] (Also noteworthy in this same time frame is the work of Ivan Ivanovich Zhegalkin (1896–1947), who, independently, provided a Boolean-valued truth-functional analysis of propositions of propositional [Zhegalkin *1927*] and its extension to first-order logic [Zhegalkin *1928-29*], undertaking to apply truth tables to the formulas of propositional calculus and first-order predicate calculus.[12])

The first instance by Peirce of a truth-functional analysis which satisfies the conditions for truth tables, but is not yet constructed in tabular form, is in his 1885 article "On the Algebra of Logic: A Contribution to the Philosophy of Notation", in which he gave a proof, using the truth table method, of what has come to be known as *Peirce Law*: $((A \rightarrow B) \rightarrow A) \rightarrow A$, his "fifth icon", whose validity he tested using truth-functional analysis. In an untitled paper written in 1902 and placed in volume 4 of the Hartshorne and Weiss edition of Peirce's *Collected Papers*, Peirce displayed the following table for three terms, *x*, *y*, *z*, writing **v** for *true* and **f** for *false* ("The Simplest Mathematics"; January 1902 ("Chapter III. The Simplest Mathematics (Logic III))", RC MS #431, January 1902; see [Peirce *1933b*, 4:260–262]).

| $x$ | $y$ | $z$ |
|---|---|---|
| **v** | **v** | **v** |
| **v** | **f** | **f** |
| **f** | **v** | **f** |
| **f** | **f** | **v** |

where $z$ is the term derived from an [undefined] logical operation on the terms $x$ and $y$.[13] In February 1909, while working on his trivalent logic, Peirce applied the tablular method to various connectives, for example negation of $x$, as $\bar{x}$, written out, in his notebook ([Peirce *1865-1909*]; see [Fisch & Turquette *1966*]), as:[14]

| $x$ | $\bar{x}$ |
|---|---|
| **V** | **F** |
| **L** | **L** |
| **F** | **V** |

where, **V**, **F**, and **L** are the truth-values true, false, and indeterminate or unknown respectively, which he called "limit".[15] Russell's rendition of Wittgenstein's tabular definition of negation, as written out on the verso of a page from Russell's transcript notes, using 'W' (wahr) and 'F' (falsch) (see [Shosky *1997*, 20]), where the negation of *p* is written out by Wittgenstein as p ◊ q, with Russell adding "= ~p", to yield: p ◊ q = ~p is

| p | q |
|---|---|
| W | W |
| W | F |
| F | W |
| F | F |

The trivalent equivalents of classical disjunction and conjunction were rendered by Peirce in that manuscript respectively as

| **Θ** | **V** | **L** | **F** |
|---|---|---|---|
| **V** | **V** | **V** | **V** |
| **L** | **V** | **L** | **L** |
| **F** | **V** | **L** | **F** |

| **Z** | **V** | **L** | **F** |
|---|---|---|---|
| **V** | **V** | **L** | **F** |
| **L** | **L** | **L** | **F** |
| **F** | **F** | **F** | **F** |

Max Fisch and Atwell R. Turquette [*1966*, 72], referring to [Turquette *1964*, 95–96], assert that the tables for trivalent logic in fact were extensions of Peirce's truth tables for bivalent logic, and hence prior to 23 February 1909 when he undertook to apply matrices for the truth-functional analysis for trivalent logic. The reference is to that part of Peirce's [*1885*, 183–193], "On the Algebra of Logic: A Contribution to the Philosophy of Notation"—§II "Non-relative Logic" —dealing with truth-functional analysis, and Turquette [*1964*, 95] uses "truth-function analysis" and "truth-table" synonymously, a confusion which, in another context, [Shosky *1997*] when warning against confusing, and insisting upon a careful distinction between the truth-table *technique* and the truth-table *device*.

Roughly contemporary with the manuscript "The Simplest Mathematics" is "Logical Tracts. No. 2. On Existential Graphs, Euler's Diagrams, and Logical Algebra", *ca*. 1903 [Peirce *1933b*, 4.476]; Harvard Lectures on Pragmatism, 1903 [Peirce *1934*, 5.108]).

In the undated manuscript [Peirce *n.d.(b)*] identified as composed *circa* 1883-84 "On the Algebra of Logic" and the accompanying supplement, we find what unequivocally would today be labeled as an indirect or abbreviated truth table for the formula {((*a* —< *b*) —< *c*) —< *d*} —< *e*, as follows:

```
{((a —< b) —< c) —< d} —< e
  f   f         f    f   —< f
  f   v        (     v )   f
  -   -   - -   -    v
```

The whole of the undated eighteen-page manuscript "Logic of Relatives", also identified as composed *circa* 1883-84 [Peirce *n.d.(c)*; MS #547], is devoted to a truth-functional analysis of the conditional, which includes the equivalent, in list form, of the truth table for *x* —< *y*, as follows [Peirce *n.d.(c)*; MS #547:16; 17]:

*x* —< *y*

| is true when | is false when |
|---|---|
| *x* = **f** *y* = **f** | *x* = **v** *y* = **f** |
| *x* = **f** *y* = **v** | |
| *x* = **v** *y* = **v** | |

Peirce also wrote follows [Peirce *n.d.(c)*; MS #547: 16] that: "It is plain that *x* —< *y* —< *z* is false only if x = v, (*y* —< *z*) = **f**, that is only if *x* = **v**, *y* = **v**, *z* = **f**…."

Finally, in the undated manuscript "An Outline Sketch of Synechistic Philosophy" identified as composed in 1893, we have an unmistakable example of a truth table matrix for a proposition and its negation [Peirce *1893*; MS #946:4], as follows:

| | t | f |
|---|---|---|
| t | t | f |
| f | t | t |

which is clearly and unmistakably equivalent to the truth-table matrix for $x \multimap y$ in the contemporary configuration, expressing the same values as we note in Peirce's list in the 1883-84 manuscript "Logic of Relatives" [Peirce *n.d.(c)*; MS #547:16; 17]. That the multiplication matrices are the most probable inspiration for Peirce's truth-table matrix is that it appears alongside matrices for a multiplicative two-term expression of linear algebra for $\{i, j\}$ and $\{i, i - j\}$ [Peirce *1893*; MS #946:4]. Indeed, it is virtually the same table, and in roughly—*i.e.*, apart from inverting the location within the respective tables for antecedent and consequent—the same configuration as that found in the notes, taken in April 1914 by Thomas Stearns Eliot (1888–1965) in Russell's Harvard University logic course (as reproduced at [Shosky *1997*, 23]), where we have:

$p \vee q$

| $q$ \ $p$ | T | F |
|---|---|---|
| T | T | T |
| F | T | F |

$p \supset q$

| $q$ \ $p$ | T | F |
|---|---|---|
| T | T | T |
| F | F | T |

$\sim p \vee \sim q$

| $q$ \ $p$ | T | F |
|---|---|---|
| T | T | T |
| F | T | F |

The ancestor of Peirce's truth table appeared thirteen years earlier, when in his lectures logic on he presented his Johns Hopkins University students with daigrammatic representations of the four combinations that two terms can take with respect to truth values. A circular array for the values $\overline{ab}$, $a\bar{b}$, $\bar{a}b$, and $ab$, each combination occupying its own quadrant:

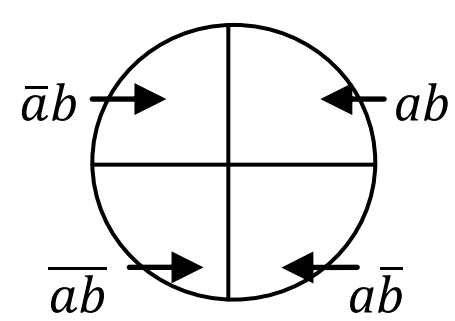

appeared in the lecture notes from the autumn of 1880 of Peirce's student Allan Marquand (1853–1924) (see editors' notes, [Peirce *1989*, 569]). An alternative array presented by Peirce himself (see editors' notes, [Peirce *1989*, 569]), and dating from the same time takes the form of a pyramid:

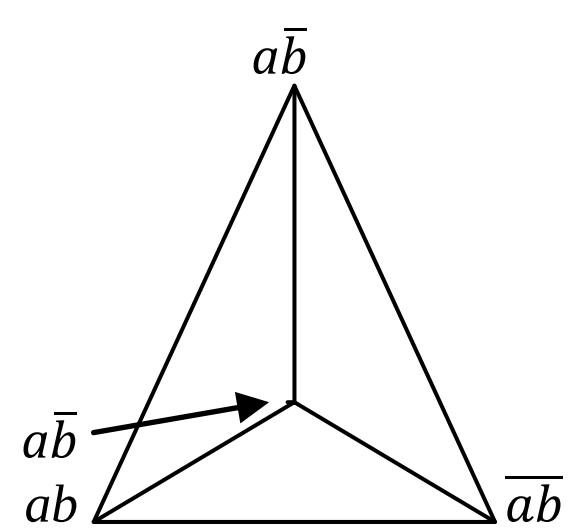

Finally, and also contemporaneous with this work, and continuing to experiment with notations, Peirce developed his "box-X" or "X-frame" notation, which resemble the square of opposition in exploring the relation between the relations between two terms or propositions. Lines may define the perimeter of the square as well as the diagonals between the vertices; the presence of connecting lines between the corners of a square or box indicates the states in which those relations are false, or invalid, absence of such a connecting line indicates relations in which true or valid relation holds. In particular, as part of this work, Peirce developed a special iconic notation for the sixteen binary connectives, as follows (from "The Simplest Mathematics" written in January 1902 ("Chapter III. The Simplest Mathematics (Logic III)", MS 431; see [Clark *1997*, 309]) containing a table presenting the 16 possible sets of truth values for a two-term proposition:

| 1 | | 2 | 3 | 4 | 5 | | 6 | 7 | 8 | 9 | 10 | 11 | | 12 | 13 | 14 | 15 | | 16 |
|---|---|---|---|---|---|---|---|---|---|---|---|---|---|---|---|---|---|---|---|
| F | | F | F | F | T | | T | T | T | F | F | F | | F | T | T | T | | T |
| F | | F | F | T | F | | T | F | F | T | T | F | | T | F | T | T | | T |
| F | | F | T | F | F | | F | T | F | T | F | T | | T | T | F | T | | T |
| F | | T | F | F | F | | F | F | T | F | T | T | | T | T | T | F | | T |

that enabled him to give a quasi-mechanical procedure for identifying thousands of tautologies from the substitution sets of expressions with up to three term (or proposition) variables and five connective variables. This table was clearly inspired by, if not actually based upon, the table presented by Christine Ladd-Franklin (*née* Ladd; (1847–1930) in her paper "On the Algebra of Logic" for the combinations $\overline{ab}, a\overline{b}, \bar{a}b,$ and $ab$ [Ladd-Franklin *1883*, 62] who, as she noted [Ladd-Franklin *1883*, 63], borrowed it, with slight modification, from Jevons's textbook, *The Principles of Science* [Jevons *1874*; *1879*, 135]. She pointed out [Ladd-Franklin *1883*, 61] that for *n* terms, there are $2^n$-many possible combinations of truth values, and she went on to provide a full-scale table for the "sixteen possible combinations of the universe with respect to two terms. Writing 0 and 1 for *false* and *true* respectively and replacing the assignment of the truth-value false with the negation of the respective terms, she arrived at her table [Ladd-Franklin *1883*, 62] providing sixteen truth values of $\{ab\}$, $\{a\overline{b}\}$, $\{\bar{a}b\}$, and $\{ab\}$.

In his X-frames notation, the open and closed quadrants are indicate truth or falsity respectively, so that for example, ⊠, the completely closed frame, represents row 1 of the table for the sixteen binary connectives, in which all assignments are false, and x, the completely open frame, represents row 16, in which all values are true (for details, see [Clark *1987*] and [Zellweger *1987*]). The X-frame notation is based on the representation of truth-values for two terms as follows:

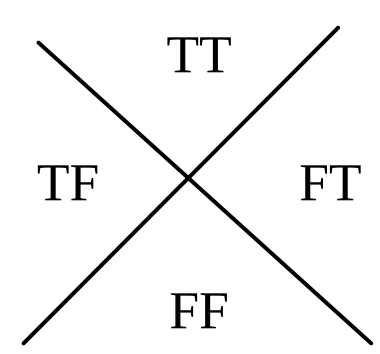

The full details of this scheme are elaborated by Peirce in his manuscript "A Proposed Logical Notation (Notation)" of *circa* 1903 [Peirce *ca. 1903*, esp. 530:26-28].

Stressing the issue of the identity of the first recognizable and ascribable example of a truth table device, or truth table matrix, the conception that it was either Wittgenstein, Post, or Łukasiewicz, or each independently of one another but almost simultaneously was challenged by Shosky [*1997*] and moved forward by a decade. But this supposition ignores the evidence advanced in behalf of Peirce going back as far as the 1930s; specifically, George W. D. Berry [*1952*] had already noted that the there is a truth table to be discovered in the work of Peirce. He was unaware of the example in Peirce's 1893 manuscript or even in the manuscripts of 1902-09 to which [Turquette 1964] and [Fisch & Turquette *1966*] point, and which were included in the Harstshorne and Weiss edition of Peirce's work. Rather, Berry referred to Peirce's published [*1885*] "On the Algebra of Logic: A Contribution to the Philosophy of Notation". There is, as we have seen, indisputablely, truth-functional analysis to be found in that work by Peirce, and an indirect truth table as well. Rather, until the inclusion of Peirce's work on trivalent logic by Hartshorne and Weiss, actual truth tables, there was no published evidence for Peirce's presentation of the truth table device. Hence, we must take Robert Lane's argument *cum grano salis* who, referring to [Berry *1952*], [Fisch & Turquette *1966*, 71–72], [Łukasiewicz & Tarski *1930*, 40, *n*. 2], and [Church *1956*, 162], in their citations of Peirce's [*1885*, 191; *1933*, 213, 4.262], when Lane [*1999*, 284] asserts that "for many years, commentators have recognized that Peirce anticipated the truth-table method for deciding whether a wff is a tautology," and agrees with Berry [*1952*, 158] that "it has long been known that [Peirce] gave an example of a two-valued truth table," explaining [Lane *1999*, 304, *n*. 4] that Berry [*1952*, 158]

> acknowledges this early appearance of the truth table. Peirce used the 1902 truth table, not to display the interpretations (or, as he himself said, the sets of values) on which a specific compound formula, consisting of three formulae, is true. He did not indicate the compound formula he had in mind. He seems to have intended the truth table to illustrate his claim that "a good many propositions concerning thee quantities cannot be expressed" using propositional connectives.

But this certainly fails to count as evidence against the claims of Fisch and Turquette to have identified truth tables in Peirce's "The Simplest Mathematics" of January 1902 and published in [Peirce *1933*, 4:260–262]. And it fails to explain why or how Shosky might have missed the evidence that "it has long been known that [Peirce] gave an example of a two-valued truth table."

What should be unconditionally recognized, in any event, is that Peirce was already well under way in devising techniques for the truth-functional analysis of propositions that these results occur quite patently and explicitly in his published work by *at least* 1885, where he was also concerned with the truth-functional analysis of the conditional, and that an unalloyed example of a truth table matrix is located in his writings that dates to at least 1893.

2. *Decomposition of propositions into function and argument instead of into subject and predicate*:

The "Booleans" were well-acquainted, from the 1820s onward with the most recent French work in function theory of their day, and although they did not explicitly employ a function-theoretical syntax in their analysis of propositions, they adopted the French algebraic approach, favored by the French analysts, Joseph-Louis Lagrange (1736–1813), Adrien-Marie Legendre (1752–1833), and Augustin-Louis Cauchy (1789–1857), to functions rather than the function-argument syntax which Frege adapted from the analysis, including in particular his teacher Karl Weierstrass (1815–1897). Moreover, Boole, De Morgan and a number of their British contemporaries who contributed to the development of "symbolical algebra" were enthusiastic adherents of this algebraic approach to analysis.[16] So there is some justification in the assertion, by Bertrand Russell, that the algebraic logicians were more concerned with logical equations than with implication (in "Recent Italian Work on the Foundations of Mathematics" of 1901; see [Russell *1993*, 353]). Ivor Grattan-Guinness [*1988*; *1997*] emphasizes the choice between algebra and function-theoretic approaches to analysis, and more generally between algebra and analysis to reinforce the distinction between algebraic logic and logistic, or function-theoretic logic, only the latter being *mathematical* logic properly so-called. This does not negate the fact, however, that algebraic logicians introduced functions into their logical calculi. If an early example is wanted, consider, e.g., Boole's definition in *An Investigation of the Laws of Thought* [Boole *1854*, 71]: "Any algebraic expression involving a symbol *x* is termed a function of *x*, and may be represented under the abbreviated general form *f*(*x*)," following which binary functions and *n*-ary functions are allowed, along with details for dealing with these as elements of logical equations in a Boolean-valued universe.

According to van Heijenoort [*1967b*, 325], Boole left his propositions unanalyzed. What he means is that propositions in Boole are mere truth-values. They are not, and *cannot be*, analyzed, until quantifiers, functions (or predicate letters), variables, and quantifiers are introduced. Even if we accept this interpretation in connection with Boole's algebraic logic, it does not apply to Peirce. We see this in the way that the Peirceans approached indexed "logical polynomials. Peirce provides quantifiers, relations, which operate as functions do for Frege, as well as variables and constants, the latter denoted by indexed terms. It is easier to understand the full implications when examined from the perspective of quantification theory. But, as preliminary, we can consider Peirce's logic of relations and how to interpret these function-theoretically.

With respect to Boole, it is correct that he conceived of propositions as adhering to the subject-predicate form and took the copula as an operator of class inclusion, differing from Aristotle only to the extent that the subject and predicate terms represented classes that were bound by no existential import, and might be empty. De Morgan, however, followed Leibniz in treating the copula as a relation rather than as representing a subsistence between an object and a property. Peirce followed De Morgan in this respect, and expanded the role of relations significantly, not merely defining the subsistence or nonsubsistence of a property in an object, but as a defined correlation between terms, such as the relation "father of" or his apparent favorite, "lover of". Boole, that is to say, followed Aristotle's emphasis on logic as a logic of classes or terms and their inclusion or noninclusion of elements of one class in another, with the copula taken as an inherence property which entailed existential import, and treated syllogisms algebraically as equations in a logic of terms. Aristotle recognized relations, but relegated them to obscurity, whereas De Morgan undertook to treat the most general form of a syllogism as a sequence of relations and their combinations, and to do so algebraically. De Morgan's algebraic logic of relations is, thus, the counterpart of Boole's algebra of classes. We may summarize the crucial distinctions by describing the core of Aristotle's formal logic as a syllogistic logic, or logic of terms and the propositions and syllogisms of the logic having a subject-predicate syntax, entirely linguistic, the principle connective for which, the copula is the copula of existence, which is metaphysically based and concerns the inherence of a property, whose reference is the predicate, in a subject; Boole's formal logic as a logic of classes, the terms of which represent classes, and the copula being the copula of class inclusion, expressed algebraically; and De Morgan's formal logic being a logic of relations whose terms are relata, the copula for which is a relation, expressed algebraically. It is possible to then say that Peirce in his development dealt with each of these logics, Aristotle's Boole's, and De Morgan's, in turn, and arrived at a formal logic which combined, and then went beyond, each of these, by allowing his copula of illation to hold, depending upon context, for terms of syllogisms, classes, and propositions, expanding these to develop (as we shall turn to in considering van Heijenoort's third condition or characteristic of the "Fregean revolution"), a quantification theory as well. Nevertheless, Gilbert Ryle (1900–1976) [*1957*, 9–10] although admittedly acknowledging that the idea of *relation* and the resulting relational inferences were "made respectable" by De Morgan, but he attributed to Russell their codification by in *The Principles of Mathematics*.

It should be borne in mind, however, that Boole did not explicitly explain how to deal with equations in terms of functions, in his *Mathematical Analysis of Logic* [Boole *1847*], although he there [Boole *1847*, 67] speaks of "elective symbols" rather than what we would today term "Boolean functions",[17] and doing so indirectly rather than explicitly. In dealing with the properies of elective functions, Boole [*1847*, 60–69] entertains Prop. 5 [Boole *1847*, 67] which, Wilfrid Hodges [*2010*, 34] calls "Boole's rule" and which, he says, is Boole's study of the deep syntactic parsing of elective symbols, and which allows us to construct an analytical tree of the structure of elective equations. Thus, for example, where Boole explains that, on considering an equation having the general form $a_1t_1 + a_2t_2 + \ldots + a_rt_r = 0$, resolvable into as many equations of the form $t = 0$ as there are non-vanishing moduli, the most general transformation of that equation is form $\psi(a_1t_1 + a_2t_2 + \ldots + a_rt_r) = \psi(0)$, provided is is taken to be of a "perfectly arbitrary character and is permitted to involve new elective symbols of any possible relation to the original elective symbols. What this entails, says Hodges [*2010*, 4] is that, given $\psi(x)$ is a Boolean function of one variable and $s$ and $t$ are Boolean terms, then we can derive $\psi(s) = \psi(s)$ from $s = t$, and, moreover, for a complex expression $\psi(x) = fghjk(x)$ are obtained by composition, such that $fghjk(x)$ is obtained by applying $f$ to $ghjk(x)$, $g$ to $hjk(x)$, …, $j$ to $k(x)$, in turn, the parsing of which yields the tree

$$\begin{array}{l} \psi(x) = f(\qquad\qquad\qquad) \\ \qquad\quad g(\qquad\qquad) \\ \qquad\quad\; h(\qquad\quad) \\ \qquad\qquad j(\qquad) \\ \qquad\qquad k(\quad) \\ \qquad\qquad\quad x \end{array}$$

in which the parsing of $\psi(s)$ and $\psi(s)$ are precisely identical except that, at the bottom node, $x$ is replaced by $s$ and $t$ respectively. If $s$ and t are also complex, then the tree will continue further. Hodges' [*2010*, 4] point is that traditional, *i.e*. Aristotelian, analysis of the syllogism makes no provision for such complexity of propositions, or, indeed, for their treatment as equations which are further analyzable beyond the simple grammar of subject and predicate.

It is also worth noting that Frege, beginning in the *Begriffsschrift* and thereafter, employed relations in a fashion similar to Peirce's. In working out his axiomatic definition of arithmetic, Frege employed the complex ancestral and proper ancestral relation to distinguish between random series of numbers from the sequence of natural numbers, utilizing the proper ancestral relation to define the latter (see [Anellis *1994*, 75–77] for a brief exposition).

The necessary apparatus to do this is provided by Ramsey's Maxim, which (in its original form), states: $x \in f \equiv f(x)$. (Recall that $f(x) = y$ is the simplest kind of mathematical expression of a *function f*, its *independent variable x*, and its *dependent variable y*, whose *value* is determined by the value of $x$. So, if $f(x) = x + 2$ and we take $x = 2$, then $y = 4$. In the expression $f(x) = y$, the *function f* takes $x$ as its *argument*, and $y$ is its *value*. Suppose that we have a binary relation $aRb$.) This is logically equivalent to the function theoretic expression $R(a, b)$, where $R$ is a binary function taking $a$ and $b$ as its arguments. A function is a relation, but a special kind of relation, then, which associates one element of the *domain* (the universe of objects or terms comprising the arguments of the function) to precisely one element of the *range*, or *codomain*, the universe of objects or terms comprising the values of the function.[18] Moreover, [Shalak *2010*] demonstrated that, for any first-order theory with equality, the domain of interpretation of which contains at least two individuals, there exists mutually embeddable theory in language with functional symbols and only one-place predicate.

In his contribution "On a New Algebra of Logic" for Peirce's [*1883a*] *Studies in Logic* of 1883 [Mitchell *1883*], his student Oscar Howard Mitchell (1851–1889) defined [Mitchell *1883*, 86] the indexed "logical polynomials", such as '$l_{i,j}$', as functions of a class of terms, in which for the logical polynomial $F$ as a function of a class of terms $a$, $b$, …, of the universe of discourse $U$, $F1$ is defined as "All $U$ is $F$" and $Fu$ is defined as "Some $U$ is $F$". Peirce defined identity in second-order logic on the basis of Leibniz's Identity of Indiscernibles, as $l_{ij}$ (meaning that every predicate is true/

false of both *i*, *j*). Peirce's quantifiers are thus distinct from Boolean connectives. They are, thus part of the "first-intensional" logic of relatives.

This takes us to the next point: that among Frege's creations that characterize what is different about the mathematical logic created by Frege and helps define the "Fregean revolution", *viz*., a quantification theory, based on a system of axioms and inference rules.

Setting aside for the moment the issue of quantification in the classical Peirce-Schröder calculus, we may summarize Peirce's contributions as consisting of a combination and unification of the logical systems of Aristotle, Boole, and De Morgan into a single system as the algebra of logic. With Aristotle, the syllogistic logic is a logic of terms; propositions are analyzed according to the subject-predicate syntactic schema; the logical connective, the copula, is the copula of existence, signalling the inherence of a property in a subject; and the syntax is based upon a linguistic approach, particularly natural language, and is founded on metaphysics. With Boole, we are presented with a logic of classes, the elements or terms of which are classes, and the logical connective or copula is class inclusion; the syntactic structure is algebraic. With De Morgan, we are given a logic of relations, whose component terms are relata; the copula is a relation, and the syntactic structure is algebraic. We find all of these elements in Peirce's logic, which he has combined within his logic of relatives. (See **Table 1**.) The fact that Peirce applied one logical connective, which he called illation, to seve as a copula holding between terms, classes, and relata, was a basis for one of the severest criticisms levelled against his logic by Bertrand Russell and others, who argued that a signal weakness of Peirce's logic was that he failed to distinguish between implication and class inclusion (see [Russell *1901c*; *1903*, 187]), referring presumably to Peirce's [*1870*],[19] while both Russell and Peano criticized Peirce's lack of distinction between class inclusion and set membership. Indeed, prior to 1885, Peirce made no distinction between sets and classes, so that Russell's criticism that Peirce failed to distinguish between class inclusion and set membership is irrelevant in any event. What Russell and Peano failed to appreciate was that Peirce intended his illation to serve as a generalized, nontransitive, copula, whose interpretation, as class inclusion, implication, or set elementhood, was determined strictly by the context in which it applied. Reciprocally, Peirce criticized Russell for failure on Russell's part to distinguish material implication and truth-functional implication (conditionality) and for his erroneous attempt to treat classes, in function-theoretic terms, as individual entities.

| **Aristotle** | **Boole** | **De Morgan** |
|---|---|---|
| syllogistic logic (logic of terms) | logic of classes | logic of relations |
| subject-predicate | classes | terms (relata) |
| copula (existence) – inherence of a property in a subject | copula (class inclusion | copula (relation) |
| linguistic/metaphysical | algebraic | algebraic |
| Peirce = Aristotle + Boole + De Morgan | | |

**Table 1. Peirce = Aristotle + Boole + De Morgan**

It is also worth noting that the development of Peirce's work followed the order from Aristotle to Boole, and then to De Morgan. That is, historically, the logic that Peirce learned was the logic of Aristotle, as it was taught at Harvard in the mid-19th century, the textbook being Richard Whately's (1787–1863) *Elements of Logic* [Whatley *1845*]. Peirce then discovered the logic of Boole, and his first efforts were an attempt, in 1867 to improve upon Boole's system, in his first publication "On an Improvement in Boole's Calculus of Logic" [Peirce *1868*]. From there he went on to study the work of De Morgan on the logic of relations and undertook to integrate De Morgan's logic of relations and Boole's algebra of logic, to develop his own logic of relatives, into which he later introduced quantifiers.[20]

3. *Peirce's quantification theory, based on a system of axioms and inference rules*:

Despite numerous historical evidences to the contrary and as suggested as long ago as by the 1950s (*e.g*.: [Berry *1952*], [Beatty *1969*], Martin *1976*]),[21] we still find, even in the very latest Peirce *Transactions*, repetition of old assertion by Quine from his *Methods of Logic* textbook [Quine *1962*, i], [Crouch *2011*, 155] that:

> In the opening sentence of his *Methods of Logic*, W. V. O. Quine writes, "Logic is an old subject, and since 1879 it has been a great one." Quine is referring to the year in which Gottlob Frege presented his *Begriffsschrift*, or "concept-script," one of the first published accounts of a logical system or calculus with quantification and a function-argument analysis of propositions. There can be no doubt as to the importance of these introductions, and, indeed, Frege's orientation and advances, if not his particular system, have proven to be highly significant for much of mathematical logic and research pertaining to the foundations of mathematics.

Quine himself ultimately acknowledged, in 1985 [Quine *1985*] and again in 1995 [Quine *1995*], that Peirce had developed a quantification theory just a few years after Frege.

Crouch is hardly alone, even at this late date and despite numerous expositions such the 1950s of Peirce's contributions, some antedating, some contemporaneous with Frege's, in maintaining the originality, and even uniqueness, of Frege's creation of mathematical logic. Thus, for example, Alexander Paul Bozzo [*2010-11*, 162] asserts and defends the historiographical phenomenon of the Fregean revolution, writing not only that Frege "is widely recognized as one of the chief progenitors of mathematical logic," but even that "Frege revolutionized the then dominant Aristotelian conception of logic," doing so by single-handedly "introducing a formal language now recognized as the predicate calculus," and explaining that: "Central to this end were Frege's insights on quantification, the notation that expressed it, the logicist program, and the extension of mathematical notions like function and argument to natural language."

It is certainly true that Peirce worked almost exclusively in equational logic until 1868.[22] But he abandoned equations after 1870 to develop quantificational logic. This effort was, however, begun as, early as 1867, and is articulated in print in "On an Improvement in Boole's Calculus of Logic [Peirce 1868]. His efforts were further enhanced by notational innovations by Mitchell in Mitchell's [*1883*] contribution to Peirce's *Studies in Logic*, "On a New Algebra of Logic", and more fully articulated and perfected, to have not only a first-, but also a second-order, quantificational theory, in Peirce's [*1885*] "On the Algebra of Logic: A Contribution to the Philosophy of Notation". Peirce himself was dissatisfied with Boole's—and others'—efforts to deal with quantifiers "some" and "all", declaring in "On the Algebra of Logic: A Contribution to the Philosophy of Notation" [Peirce *1885*, 194] that, until he and Mitchell devised their notation in 1883, no one was able to properly handle quantifiers, that: All attempts to introduce this distinction into the Boolian algebra were more or less complete failures until Mr. Mitchell showed how it was to be effected."

But, even more importantly, that Peirce's system dominated logic in the final two decades of the $19^{th}$ century and first two decades of the $20^{th}$.

By 1885, Peirce not only had a fully developed first-order theory, which he called the *icon of the second intention*, but a good beginning at a second-order theory. Our source here is Peirce's [*1885*] "On the Algebra of Logic: A Contribution to the Philosophy of Notation". In "Second Intentional Logic" of 1893 (see [Peirce *1933b*, 4.56–58), Peirce even presented a fully developed second-order theory.

The final version of Peirce's first-order theory uses indices for enumerating and distinguishing the objects considered in the Boolean part of an equation as well as indices for quantifiers, a concept taken from Mitchell.

Peirce introduced indexed quantifiers in "The Logic of Relatives" [Peirce *1883b*, 189]. He denoted the existential and universal quantifiers by '$\Sigma_i$' and '$\Pi_i$' respectively, as logical sums and products, and individual variables, $i$, $j$, …, are assigned both to quantifiers and predicates. He then wrote '$l_{i,j}$' for '*i* is the lover of *j*'. Then "Everybody loves somebody" is written in Peirce's quantified logic of relatives as $\Pi_i \, \Sigma_j \, l_{i,j}$, i.e. as "Everybody is the lover of somebody". In Peirce's own exact expression, as found in his "On the Logic of Relatives" [*1883b*, 200]), we have: "$\Pi_i \, \Sigma_j \, l_{i,j} > 0$ means that everything is a lover of something." Peirce's introduction of indexed quantifier in fact establishes Peirce's quantification theory as a many-sorted logic.

That is, Peirce defined the existential and universal quantifiers, in his mature work, by '$\Sigma_i$' and '$\Pi_i$' respectively, as logical sums and products, e.g., $\Sigma_i x_i = x_i + x_j + x_k + \ldots$, and $\Pi_i x_i = x_i \bullet x_j \bullet x_k$, and individual variables, $i$, $j$, …, are assigned both to quantifiers and predicates. (In the Peano-Russell notation, these are $(\exists x)F(x) = F(x_i) \vee F(x_j) \vee F(x_k)$ and are $(\forall x)F(x) = F(x_i) \bullet F(x_j) \bullet F(x_k)$ respectively.)

The difference between the Peirce-Mitchell-Schröder formulation, then, of quantified propositions, is purely cosmetic, and both are significantly notationally simpler than Frege's. Frege's rendition of the proposition "For all $x$, if $x$ is $F$, then $x$ is $G$", i.e. $(\forall x)[F(x) \supset G(x)]$, for example, is and, in the Peirce-Mitchell-Schröder notation could be formulated as $\Pi_i \, (f_i \prec g_i)$, while "There exists an $x$ such that $x$ is $f$ and $x$ is $G$, in the familiar Peano-Russell notation is formulated as $(\exists x)[F(x) \bullet G(x)]$, and $\Sigma_i \, (f_i \prec g_i)$ in the Peirce-Mitchell-Schröder notation, is rendered as in Frege's notation, that is ~(∀*x*)~[*F*(*x*) • *G*(*x*)].

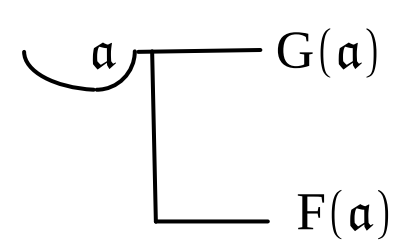

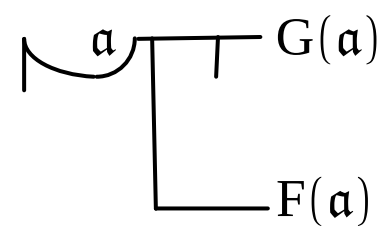

Not only that; recently, Calixto Badesa [*1991*; *2004*], Geraldine Brady [*2000*] (see also [Anellis *2004b*]) in detail, and Enrique Casanovas [*2000*] briefly and emphasizing more recent developments, traced the development of the origins of the special branches of modern mathematical logic known as *model theory*,[23] which is concerned with the properties of the consistency, completeness, and independence of mathematical theories, including of course the various logical systems, and *proof theory*, concerned with studying the soundness of proofs within a mathematical or logical system. This route runs from Peirce and his student Oscar Howard Mitchell (1851–1889) through Ernst Schröder (1841–1902) to Leopold Löwenheim (1878–1957), in his [*1915*] "Über Möglichkeiten im Relativkalkul", Thoralf Skolem (1887–1963), especially in his [*1923*] "Einige Bemerkungen zur axiomatischen Begründung der Mengenlehre", and—I would add—Jacques Herbrand (1908–1931), in his [*1930*] *Recherches sur la théorie des démonstration*. It was based upon the Peirce-Mitchell technique for elimination of quantifiers by quantifier expansion that the Löwenheim-Skolem Theorem (hereafter LST) allows logicians to determine the validity of within a theory of the formulas of the theory and is in turn the basis for Herbrand's Fundamental Theorem (hereafter FT), which can best be understood as a strong version of LST. In his article "Logic in the Twenties:" Warren D. Goldfarb recognized that Peirce and Schröder contributed to the development of quantification theory and thus states that [Goldfarb *1979*, 354]:

> Building on earlier work of Peirce, in the third volume of his **Lectures on the algebra of logic** [1895] Schröder develops the calculus of relatives (that is, relations). Quantifiers are defined as certain possibly infinite sums and products, over individuals or over relations. There is no notion of formal proof. Rather, the following sort of question is investigated: given an equation between two expressions of the calculus, can that equation be satisfied in various domains—that is, are there relations on the domain that make the equation true? This is like our notion of satisfiability of logical formulas. …Schröder seeks to put the entire algebra into a form involving just set-theoretic operations on relations (relative product and the like), and no use of quantifiers. In this form, the relation-signs stand alone, with no subscripts for their arguments, so that the connection between the calculus and predication tends to disappear. (In order to eliminate subscripts Schröder often simply treats individuals as though they were relations in their own right. This tactic leads to confusions and mistakes; see Frege's review [1895].) Schröder's concern is with the laws about relations on various universes, and not with the expressive power he gains by defining quantification in his (admittedly shady) algebraic manner.

(Note first of all Goldfarb's blatantly dismissive manner, and his abusive—"shady"—treatment of Schröder. But notice also that his concern is preeminently with what he conceives as Schröder's understanding of the foremost purpose of the algebra of logic as a calculus, and that Peirce disappears immediately from Goldfarb's ken. But what is relevant for us in Goldfarb's [*1979*] consideration, which is his aim, namely to establish that the expansion of quantified equations into their Boolean equivalents, that is, as sums and products of a finite universe prepared the way for the work of Löwenheim and Skolem. Not dismissive, but rather oblivious, to the historical roots of the Löwenheim-Skolem Theorem, Herbrand's Fundamental Theorem, and related results in proof theory and model theory despite its strong concern for the historical background of these areas of logic, [Hoffmann *2011*], despite a reference to [Brady *2000*], fails to mention algebraic logic at all in its historical discussion, treating the history of the subject exclusively from the perspective of the Fregean position and mentioning Schröder only in connection with the Cantor-Schröder-Bernstein Theorem of set theory, and Peirce not at all.)

The original version of what came to be known as the LST, as stated by Löwenheim, is simply that:

*If a well-formed formula of first-order predicate logic is satisfiable, then it is* $\aleph_0$*-satisfiable.*

Not only that: in the manuscript "The Logic of Relatives: Qualitative and Quantitative" of 1886, Peirce himself is making use of what is essentially a finite version of LST,[24] that

*If F is satisfiable in every domain, then F is* $\aleph_0$*-satisfiable*;

that is:

*If F is n-satisfiable, then F is (n* + 1*)-satisfiable,*

and indeed, his proof was in all respects similar to that which appeared in Löwenheim's 1915 paper, where, for any $\kappa < \lambda$, a product vanishes (*i.e*. is satisfiable), its $\kappa^{th}$ term vanishes.

In its most modern and strictest form, the LST says that:

*For a* $\boldsymbol{\kappa}$*-ary universe, a well-formed formula F is* $\aleph_0$*-valid if it is* $\kappa$*-valid for every finite* $\boldsymbol{\kappa}$*, provided there is no finite domain in which it is invalid*.

Herbrand's FT was developed in order to answer the question: what finite sense can generally be ascribed to the truth property of a formula with quantifiers, particularly the existential quantifier, in an infinite universe? The modern statement of FT is:

*For some formula F of classical quantification theory, an infinite sequence of quantifier-free formulas* $F_1$*,* $F_2$*, ... can be effectively generated, for F provable in (any standard) quantification theory, if and only if there exists a* $\kappa$ *such that F is (sententially) valid; and a proof of F can be obtained from* $F_\kappa$.

For the role that Peirce's formulation of quantifier theory played in the work of Schröder and the algebraic logicians who followed him, as well as the impact which it had more widely, not only on Löwenheim, Skolem, and Herbrand, Hilary Putnam [*1982*, 297] therefore conceded at least that: "Frege did "discover" the quantifier in the sense of having the rightful claim to priority. But Peirce and his students discovered it in the effective sense."

The classical Peirce-Schröder calculus also, of course, played a significant role in furthering the developments in the twentieth century of algebraic itself, led by Alfred Tarski (1903–1983) and his students. Tarski was fully cognizant from the outset of the significance of the work of Peirce and Schröder, writing, for example in "On the Calculus of Relations" [Tarski *1941*, 73–74], that:

> The title of creator of the theory of relations was reserved for C. S. Peirce. In several papers published between 1870 and 1882, he introduced and made precise all the fundamental concepts of the theory of relations and formulated and established its fundamental laws. Thus Peirce laid the foundation for the theory of relations as a deductive discipline; moreover he initiated the discussion of more profound problems in this domain. In particular, his investigations made it clear that a large part of the theory of relations can be presented as a calculus which is formally much like the calculus of classes developed by G. Boole and W. S. Jevons, but which greatly exceeds it in richness of expression and is therefore incomparably more interesting from the deductive point of view. Peirce's work was continued and extended in a very thorough and systematic way by E. Schröder. The latter's *Algebra und Logik der Relative*, which appeared in 1895 as the third volume of his *Vorlesungen über die Algebra der Logik*, is so far the only exhaustive account of the calculus of relations. At the same time, this book contains a wealth of unsolved problems, and seems to indicate the direction for future investigations.
>
> It is therefore rather amazing that Peirce and Schröder did not have many followers. It is true that A. N. Whitehead and B. Russell, in *Principia mathematica*, included the theory of relations in the whole of logic, made this theory a central part of their logical system, and introduced many new and important concepts connected with the concept of relation. Most of these concepts do not belong, however, to the theory of relations proper but rather establish relations between this theory and other parts of logic: *Principia mathematica* contributed but slightly to the intrinsic development of the theory of relations as an independent deductive discipline. In general, it must be said that—though the significance of the theory of relations is universally recognized today—this theory, especially the calculus of relations, is now in practically the same stage of development as that in which it was forty-five years ago.

Thus, as we see, Tarski credited Peirce with the invention of the calculus of binary relations. Although we might trace the bare beginnings to De Morgan, and in particular to the fourth installment of his On the Syllogism" [De Morgan *1860*], Tarski [*1941*, 73] held that De Morgan nevertheless "cannot be regarded as the creator of the modern theory of relations, since he did not possess an adequate apparatus for treating the subject in which he was interested, and was apparently unable to create such an apparatus. His investigations on relations show a lack of clarity and rigor which perhaps ac-counts for the neglect into which they fell in the following years. The title of creator of the theory of relations was reserved for." In any case, it was, as Vaughn Pratt [*1992*] correctly noted, Peirce who, taking up the subject directly from De Morgan in "Description of a Notation for the Logic of Relatives…" [Peirce *1870*], brought to light the full power of that calculus, by articulating those technicalities at which De Morgan's work only hinted. More to the point for our purposes, and in contradistinction to the judgment of van Heijenoort on behalf of the "Fregeans" against the "Booleans", Pratt [*1992*, 248; my emphasis] stresses the point, among others, that: "The calculus of binary relations as we understand it today…*is a logic*."

The algebraic logic and logic of relations with which we are familiar today is the work largely of Tarski and his students, initiated by Tarski in picking up where Peirce and Schröder left off. Major results established by Tarski and his student Steven Givant in their *A Formalization of Set Theory without Variables* [Tarski & Givant *1987*] was to answer a question posed by Schröder and answered, in the negative by Alwin Reinhold Korselt (1864–1947) (but reported on Korselt's behalf by Löwenheim [*1914*, 448], as to whether every formula of binary first-order quantificational logic is reducible (expressible) in the Peirce-Schröder calculus, as stated in the Korselt-Tarski Theorem; a generalization, which asserts that set theory cannot be formulated within a theory having three two-place predicate constants of the second type and four two-place function constants of the third type, or any extensions of such a language; and the special, closely related result, the so-called Main Mapping Theorem, asserting that there is a formula of first-order quantification theory having four constants, which cannot be expressed in the three-constant theory or any of its extensions, thus—apparently—specifically contradicting Peirce's Reduction Thesis, that every equation of the logic of relations in which there is a quaternary relation can be expressed by an equation composed of a combination of monadic, dyadic, and triadic relations, by exhibiting an equation having four or more terms is reducible to an expression comprised of some combination of statements of one, two, and three terms (see [Anellis *1997*]; see [Anellis & Houser *1991*], [Maddux *1991*] and [Moore *1992*] for general historical background).

For the role that Peirce's formulation of quantifier theory played in the work of Schröder and the algebraic logicians who followed him, as well as the impact which it had more widely, not only on Löwenheim, Skolem, and Herbrand, Hilary Putnam [*1982*, 297] therefore remarked that: "Frege did "discover" the quantifier in the sense of having the rightful claim to priority. But Peirce and his students discovered it in the effective sense."

4. *Peirce's definitions of* infinite sequence *and* natural number *in terms of logical notions (i.e. the logicization of mathematics)*:

In his "The Logic of Number" [Peirce *1881*] of 1881, Peirce set forth an axiomatization of number theory, starting from his definition of *finite set* to obtain natural numbers. Given a set *N* and *R* a relation on *N*, with 1 an element of *N*; with definitions of *minimum*, *maximum*, and *predecessor* with respect to *R* and *N* given, Peirce's axioms, in modern terminology, are:

1. *N* is partially ordered by *R*.
2. *N* is connected by *R*.
3. *N* is closed with respect to predecessors.
4. 1 is the minimum element in *N*; *N* has no maximum.
5. Mathematical induction holds for *N*.

It is in this context important to consider Sluga's testimony [Sluga *1980*, 96–128], that it took five years beyond the completion date of December 18, 1878 for the *Begriffsschrift* to provide the promised elucidation of the concept of *number* following his recognition that are logical objects and realizing that he had not successfully incorporated that recognition into the *Begriffsschrift* [Frege *1879*]. Certainly, if Peirce in 1881 had not yet developed a complete and coherent logical theory of number, neither, then, had Frege before 1884 in the *Die Grundlagen der Arithmetik* [Frege *1884*]. Thus, by the same token and rationale that we would contest whether Peirce had devised a number theory within his logical system by 1881, so would we have to coext whether Frege had done so by 1879.

We have already noted, in consideration of the second criterion for describing the Fregean innovations, that Frege employed the relations of ancestral and proper ancestral, akin to Peirce's "father of" and "lover of" relation. Indeed, the Fregean *ancestrals* can best be understood as the inverse of the *predecessor* relation (for an account of the logical construction by Frege of his ancestral and proper ancestral, see [Anellis *1994*, 75–77]). Frege used his definition of the proper ancestral to define mathematical induction. If Peirce's treatment of predecessor and induction is inadequate in "The Logic of Number", it is necessary to likewise note that Frege's definition "*Nx*" of "*x* is a natural number" does not yet appear in the *Begriffsschrift* of 1879 although the definition of proper ancestral is already present in §26 of the *Begriffsschrift*. Frege's definition of "*Nx*", as van Heijenoort himself remarked, does not appear until Part II of the *Grundlagen der Arithmetik* [Frege *1884*].

The only significant differences between axiomatization of number theory by Richard Dedekind (1831–1914) and Peirce's was that Dedekind, in his [*1888*] *Was sind und was sollen die Zahlen?* started from *infinite sets* rather than finite sets in defining natural numbers, and that Dedekind is explicitly and specifically concerned with the real number continuum, that is, with infinite sets. Peirce also questioned Cantor's (and Dedekind's) approach to set theory, and in particular Cantor's conception of infinity and of the continuum, arguing for example, in favor of infinitesimals rather than the possibility of the existence of an infinite number of irrationals between the rational numbers, arguing that there cannot be as many real numbers as there are points on an infinite line, and, more importantly, arguing that Cantor's arguments in defining the real continuum as transfinite are mathematical rather than grounded in logic. Affirming his admiration for Cantor and his work nevertheless, from at least 1893 forward, Peirce undertook his own arguments and construction of the continuum, albeit in a sporadic and unsystematic fashion. Some of his results, including material dating from *circa* 1895 and from May 1908, unpublished and incomplete, were only very recently published, namely "The Logic of Quantity", and, from 1908, the "Addition" and "Supplement" to his 1897 "The Logic of Relatives" (esp. [Peirce *1897*, 206]; see, *e.g*. [Peirce *2010*, 108–112; 218–219, 221–225]). He speaks, for example of Cantor's conception of the continuum as identical with his own pseudo-continuum, writing [Peirce *2010*, 209] that "I define a *pseudo-continuum* as that which modern writers on the theory of functions call a *continuum*. But this is fully represented by, and according to G. Cantor stands in one to one correspondence with the totality [of] real values, rational and irrational…."

The equivalence of Peirce's axiomatization of the natural numbers to Dedekind's and Peano's is demonstrated by Paul Bartram Shields in his doctoral thesis *Charles S. Peirce on the Logic of Number* [Shields *1981*] and his [*1997*] paper based on the thesis.

Nevertheless, Francesco Gana [*1985*] examined the claim by Peirce that Dedekind in his *Was sind und was sollen die Zahlen?* plagiarized his "Logic of Number", and concludes that the charge was unjustified.

The technical centerpiece of Dedekind's mathematical work was in number theory especially algebraic number theory. His primary motivation was to provide a foundation for mathematics and in particular to find a rigorous definition of the real numbers and of the number continuum upon which to establish analysis in the style of Weierstrass. This means that he sought to axiomatize the theory of numbers, based upon a rigorous definition of the real number system which could be employed in defining the theory of limits of a function for use in the differential and integral calculus, real analysis, and related areas of function theory. His concern, in short, is with the rigorization and arithmetization of analysis.

For Peirce, the object behind his axiomatization of the system of natural numbers was stated in "On the Logic of Number" [Peirce *1881*, 85] as establishing that "the elementary propositions concerning number…are rendered [true] by the usual demonstrations." He therefore undertakes "to show that they are strictly syllogistic consequences from a few primary propositions," and he asserts "the logical origin of these latter, which I here regard as definitions," but for the time being takes as given. In short, Peirce here wants to establish that the system of natural numbers can be developed axiomatically by deductive methods (*i.e*. "syllogistically", applying his logic of relations, and that the system of natural numbers can be constructed by this means on the basis of some logical definitions. Peirce's central concern in "The Logic of Number" was with natural numbers and foundations of arithmetic, rather than analysis, and after he studied Dedekind's work and that of Georg Cantor (1845–1918), he began to focus deeper attention on infinite sets, but did not publish his work[25] (The comparison between Peirce's and Dedekind's approaches to set theory as evidenced in their respective axiomatizations of number, may be summarized in **Table 2**.)

Whether this is tantamount, from the philosopher's standpoint, to the logicism of Frege and Peano has been the subject of debate.[26] Peirce's declaration that mathematics is not subject to logic, whereas logic depends upon mathematics, and that the only concern that logic has with mathematical reasoning is to describe it (see [Peirce *2010*, 24; 33; *1931*, 1:191-192]), strongly suggests that Peirce was not at all a philosophical logicist in any sense that would be recognized by Dedekind, Frege or Russell or any like-minded philosophers of mathematics. Indeed, Peirce [*2010*, 32] mentions his discussions with Dedekind on their differing conceptions of the nature of the relation between logic and mathematics, and in particular his disagreement with Dedekind's position that mathematics is a branch of logic, whereas Peirce upheld the position that mathematical, or formal, logic, which is only a small portion of the whole subject of logic, and not even its principal part, is best compared with such applications of mathematics to fields outside of mathematics, as an *aid* to these fields, such as mathematical optics or mathematical economics respectively, which, after all, remains optics and economics, rather than thereby being, or becoming, mathematics [Peirce *2010*, 32–33; *1931*, 1:191–192].

As a notational aside, it may be worth recollecting here that, just as Löwenheim, Skolem, and Herbrand adopted the Peirce-Mitchell-Schröder treatment of quantifiers in their work, Ernst Friedrich Ferdinand Zermelo (1871–1953), in his full axiomatization of set theory [Zermelo *1908*] used Schröder's subsumption as the basic relation for subsets. This is of particular interest because, for Zermelo, the rationale for his effort was the Peircean one of developing an axiomatic system for number theory.

| **Peirce** | **Dedekind** |
|---|---|
| motivation (as stated on "On the Logic of Number" [Peirce *1880*, 85]: establish that "the elementary propositions concerning number…are rendered [true] by the usual demonstrations" | motivation (*Was sind und was sollen die Zahlen?* [Dedekind *1888*]: foundations for mathematics |
| natural numbers | number theory, especially algebraic numbers |
| "show that numbers are strictly syllogistic consequences from a few primary propositions" | give rigorous definition of reals & of number continuum to establish analysis in the style of Weierstrass |
| "the logical origin" of these propositions" are regarded as definitions | use this def. to define the theory of limits of a function for use in differential & integral calculus; real analysis; related areas of function therory; *i.e*., rigorize & arithmeticize analysis |
| start from *finite* sets in defining natural numbers | start from *infinite sets* in defining natural numbers |
| explicitly and specifically concerned with natural numbers and arithmetic (only later dealing with transfinite set theory, but not publishing) | explicitly and specifically concerned with the real number continuum, that is, with infinite sets |

**Table 2. Peirce's *vs.* Dedekind's Axiomatization of Number**

5. *Presentation and clarification of the concept of* formal system:

Did Peirce formally and explicitly set forth his conception of a formal system?

I would suggest that, even if Peirce nowhere formally and explicitly set forth his conception of a formal system, it is present and implicit in much of his work, in "On the Logic of Number" for example, in the explication of the purpose of his project of deducing, in a logically coherent manner, and in strict accordance with deductive inference rules on the basis of a few essential and carefully chosen and well-defined "primary propositions"—definitions, the propositions requisite for deriving and expressing the elementary propositions—axioms—of mathematics concerning numbers.

Taking a cue from Geraldine Brady, who wrote [Brady *2000*, 14] of Peirce's "failure to provide a formal system for logic, in the sense of Frege's. The motivation to create a formal system is lacking in Peirce…," and he thus "made no early attempt at an all-encompassing formal system," we must admittedly note that here is in Peirce no one set of axioms by which to derive all of logic, still less, all of mathematics. Rather, what we have is an on-going experiment in developing the basics of a logic that answers to specific purposes and has as its ultimate goal the creation of a calculus that serves as a tool for the wider conception of logic as a theory of signs. Nevertheless, this absence of a single axiomatic system designed to encompass and express all of mathematics, not only in Peirce, but in the "Boolean" tradition, will become the basis for van Heijenoort and others to argue that there is a lack of universality in the algebraic logic tradition of Boole, Peirce, and Schröder, and, consequently, that Frege was correct in asserting (*e.g.* at [Frege *1880/81*; *1882*; *1883*] that the Booleans have developed logic as a calculus, but not logic as a language. Brady's judgment in regard to Peirce's failure to develop a formal system is shared by Jaakko Hintikka [*1997*], who relies upon van Heijenoort's distinction between logic as calculus (or *logica utens*) and logic as language (or *logica docens*), John Brynes [*1998*], and others.

What we can say, I would suggest, is that Peirce undertook through the length of his career to develop a series of formal systems, without, however, expressing in detail or specifically, the concept of a formal system, certainly not in the language that would sound familiar to readers of Frege, Hilbert, Peano, or Russell. That Peirce would experiment with constructing various formal systems, without attempting to devise one formal system that was undertaken to do general duty in the same way that Frege's *Begriffsschrift* or Whitehead and Russell's *Principia Mathematica* were intended to do, may be explained in terms of the differences between a *logica utens* and a *logica docens*, and Peirce's broader conception of logic within the architectonic within which he placed the various sciences. For Peirce, as for his fellow "Booleans", it would seem that his chief concerns were for devising a *logica docens* or series of such logics. (I leave the Peirce's architectonic to philosophers to discuss. Peirce's conception of *logica utens* and *logica docens*, however, is a consideration in the discussion of the next of van Heijenoort's characteristics defining the Fregean revolution.)

6. *Making possible, and giving, a use of logic for philosophical investigations (especially for philosophy of language)*:

Van Heijenoort's understanding of Frege's conception of application of his logical theory for philosophical investigations and in particular for philosophy of language can be seen as two-fold, although van Heijenoort in particular instances envisioned it in terms of analytic philosophy. On the one hand, Frege's logicist program was understood as the centerpiece, and concerned the articulation of sciences, mathematics included, developed within the structure of the logical theory; on the other hand, it is understood, more broadly, as developing the logical theory as a universal language.

Distinguishing *logic as calculus* and *logic as language*, the "Booleans" or algebraic logicians are understood to treat logic as a mere calculus, whereas Frege and the "Fregeans" see their logic to be both a calculus and a language, but first and foremost as a language. It is in this regard that Frege criticized Schröder, although he had the entire algebraic logic tradition in mind, from Boole to Schröder. This was in response to Schröder's assertion, in his review of 1880 of Frege's *Begriffsschrift*, that Frege's system "does not differ essentially from Boole's formula language," adding: "With regard to its major content the *Begriffsschrift* could actually be considered a *transcription*"—"Umschreibung"—"of the Boolean formula language."

Peirce used logic in several different senses. In the narrowest sense, it refers to deductive logic, and is essentially equivalent to the calculus of the logic of relations. In the broader sense, it is virtually synonymous with the *grammatica speculativa*, which includes three branches: semiotics, rhetoric, and logic. In a related use, he defined logic in the broader sense as coextensive with semiotics, the theory of signs, writing [Peirce *1932*, 2.92] that: "Logic is the science of general necessary laws of signs" and [Peirce *1932*, 2. 227] that: "Logic, in its general sense, is…only another name for semiotic, the quasi-necessary, or formal, doctrine of signs."

In the narrow sense, logic is a *normative* science, establishing the rules for correctly drawing, or *deducing*, conclusions from given propositions. It is on this basis that Peirce was able, as we have seen, to translate the Aristotelian syllogism as an implication. Thus: "To *draw necessary* conclusions is one thing, to draw conclusions is another, and the *science* of drawing conclusions is another; and that science is Logic." Logic in this usage is a deductive methodology,[27] and in that case a system of logical symbols is the means by which we can "analyze a reasoning into its last elementary steps" [Peirce *1933b*, 4.239]. In an unpublished manuscript on "Logic as the Study of Signs" of 1873, intended as part of a larger work on logic, Peirce went so far as to defined logic as a study of signs. He then wrote (see "Of Logic as the Study of Signs"; MS 221; Robin catalog # 380; March 14, 1873; published: [Peirce *1986*, 82–84]). In that work Peirce explores the nature of logic as algebra, or *critic* (see [Bergman & Paavola *2003*-], "Critic, Speculative Critic, Logical Critic") and its relation with the broader field of semiotics, or *grammatica speculativa*. He then writes [Peirce *1986*, 84]:

> The business of Algebra in its most general signification is to exhibit the manner of tracing the consequences of supposing that certain signs are subject to certain laws. And it is therefore to be regarded as a part of Logic. Algebraic symbols have been made use of by all logicians from the time of Aristotle, and probably earlier. Of late, certain logicians of some popular repute, but who represent less than any other school the logic of modern science, have objected that Algebra is exclusively the science of quantity, and is therefore entirely inapplicable to Logic. This argument is not so weak that I am astonished at these writers making use of it, but it is open to three objections: In the first place, Algebra is not a science of quantity exclusively, as every mathematician knows; in the second place these writers themselves hold that logic is a science of quantity; and in the third place, they, themselves, make a very copious use of algebraic symbols in Logic.

[Anellis *forthcoming (a)*] is an attempt at an explanation of the relation of Peirce's view with that of van Heijenoort regarding logic as calculus and logic as language, in an effort to understand whether, and if so, how, Peirce may have contributed to the conception of the role of logic as a language as well as of logic as a calculus, and along the way whether logic can therefore satisfy, to some extent or not, the place of his logic in philosophically or logico-linguistic investigations; [Anellis *2011*] examines Peirce's conception of the relations between logic and language, in particular against the background of the views and attitudes specifically of Russell's contemporaries, and from the perspective of van Heijenoort's distinction logic as calculus/logic as language distinction.

We may think of logic as calculus and logic as language in terms, borrowed from the medievals of *logica utens* and *logica docens*.

In the terms formulated by van Heijenoort (see, *e.g*. [van Heijenoort *1967b*]), a *logica utens* operates with a specific, narrowly defined and fixed universe of discourse, and consequently serves as a *logic as calculus*, and thus as a *calculus ratiocinator*, whereas a *logica docens* operates with a universal domain, or universal universe of discourse, characterized by Frege as the *Universum*, which is in fact universal and fixed.[28]

For Peirce, the distinction between *logica docens* and *logica utens* was consistently formulated in terms of the *logica utens* as a "logical theory" or "logical doctrine" as a means for determining between good and bad reasoning (see, e.g. "The Proper Treatment of Hypotheses: a Preliminary Chapter, toward an Examination of Hume's Argument against Miracles, in its Logic and in its History" (MS 692, 1901) [Peirce *1932*, 2:891–892]; from the "Minute Logic", "General and Historical Survey of Logic. Why Study Logic? Logica Utens", *ca*. 1902 [Peirce *1932*, 2.186]; "Logical Tracts. No. 2. On Existential Graphs, Euler's Diagrams, and Logical Algebra", *ca*. 1903 [Peirce *1933b*, 4.476]; Harvard Lectures on Pragmatism, 1903 [Peirce *1934*, 5.108]), and the *logica docens* in terms of specific cases. In the entry on "Logic" for Baldwin's *Dictionary* [Peirce & Ladd-Franklin *1902*, II, 21], Peirce, in collaboration with his former student Christine Ladd-Franklin (1847–1930), wrote:

> In all reasoning, therefore, there is a more or less conscious reference to a general method, implying some commencement of such a classification of arguments as the logician attempts. Such a classification of arguments, antecedent to any systematic study of the subject, is called the reasoner's *logica utens*, in contradistinction to the result of the scientific study, which is called *logica docens.* See REASONING.
>
> That part of logic, that is, of *logica docens*, which, setting out with such assumptions as that every assertion is either true or false, and not both, and that some propositions may be recognized to be true, studies the constituent parts of arguments and produces a classification of arguments such as is above described, is often considered to embrace the whole of logic; but a more correct designation is Critic (Gr. κριτική. According to Diogenes Laertius, Aristotle divided logic into three parts, of which one was πρὸς κρίσιν). …

In the next paragraph, Peirce and Ladd-Franklin establish the connection between logic as critic and the *grammatica speculativa*:

> It is generally admitted that there is a doctrine which properly antecedes what we have called critic. It considers, for example, in what sense and how there can be any true proposition and false proposition, and what are the general conditions to which thought or signs of any kind must conform in order to assert anything. Kant, who first raised these questions to prominence, called this doctrine transcendentale Elementarlehre, and made it a large part of his *Critic of the Pure Reason.* But the *Grammatica Speculativa* of Scotus is an earlier and interesting attempt. The common German word is Erkenntnisstheorie, sometimes translated EPISTEMOLOGY (q.v.).

Ahti-Veikko Pietarinen [Pietarinen *2005*] has characterized the distinction for Peirce as one between the *logica utens* as a logic of action or use and the *logica docens* as a general theory of correct reasoning. In the terms formulated by Jean van Heijenoort (1912–1986) (see, *e.g*. [van Heijenoort *1967a*]), a *logica utens* operates with a specific, narrowly defined and fixed universe of discourse, and consequently serves as a *logic as calculus*, and thus as a *calculus ratiocinator*, whereas a *logica docens* operates with a universal domain, or universal universe of discourse, characterized by Frege as the *Universum*, which is in fact universal and fixed. There are several interlocking layers to van Heijenoort's thesis that, as a result of its universality, it is not possible to raise or deal with metalogical, *i.e*. metasystematic, properties of the logical system of *Principia Mathematica*. These aspects were dealt with in a series of papers by van Heijenoort over the course of more than a decade. The writings in question, among the most relevant, include "Logic as calculus and Logic as Language" [van Heijenoort *1967a*], "Historical Development of Modern Logic (1974) [van Heijenoort *1992*]; "Set-theoretic Semantics", [van Heijenoort *1977*], "Absolutism and Relativism in Logic" (1979) [van Heijenoort *1986*], and "Système et métasystème chez Russell" [van Heijenoort *1987*].

At the same time, however, we are obliged to recognize that Peirce's own understanding of *logica utens* and *logica docens* is not precisely the same as we have represented them here and as understood by van Heijenoort. For Peirce, a *logica utens* is, or corresponds to a logical theory, or logic as critic, and *logica docens* is the result of the scientific study, and more akin to an uncritically held but deeply effective logical theory, and hence normative, which because it governs our actions almost instinctively amounts almost to a moral theory.

We should distinguish more carefully *logica docens* from *logica utens* as conceived by van Heijenoort as it relates to Peirce.[29] A *logica utens* is specific calculus designed to serve a specific purpose or narrow field of operation, and is typically associated to one *universe of discourse* (a term coined by De Morgan) which applies to a specific, well-defined domain and which Schröder came subsequently to conceive as a *Denkbereich*.[30] The classical Boole-Schröder algebra of logic is understood by van Heijenoort as a *logica utens* in this sense. Although admittedly the universe of discourse can have more than one semantic interpretation, that interpretation is decided ahead of time, to apply specifically to sets, or to classes, or to propositions, but never does duty for more than one of these at a time. In a more practical sense, we might consider the axiomatic systems developed by the postulate theorists, who set forth a specific system of axioms for specific fields of mathematics, and for whom, accordingly, the universe of discourse is circumspect in accordance with the area of mathematics for which an axiomatic system was formulated. For example, we see one for group theory, and another for geometry; even more narrowly, we find one for metric geometry, another for descriptive geometry; etc. The universe of discourse for the appropriate postulate system (or *logica utens*) for geometry would consist of points, lines, and planes; another universe of discourse, might, to borrow Hilbert's famous example be populated by tables, chairs, and beer mugs. We may, correspondingly, understand the *logica docens* as an all-purpose logical calculus which does not, therefore, operate with one and only one or narrowly constrained specific universe of discourse or small group of distinct universes of discourse. Boole's and De Morgan's respective calculi, as well as Peirce's and Schröder's are *extensional*, their syntactical components being comprised of classes, although, in Peirce's logic of relations, capable likewise of being comprised of sets.

For van Heijenoort, "semantics", Juliet Floyd [*1998*, 143] insists, means either model-theoretic or set-theoretic semantics. To make sense of this assertion, we need to understand this as saying that the interpretation of the syntax of the logical language depends upon a universe of discourse, an extensional rather than intensional universe of discourse, in this case the universal domain, or *Universum*, satisfying either Frege's *Werthverlauf*, or course-of-values semantic and Russell's set-theoretic semantic. The model-theoretic approach makes no extra-systematic assumptions and is entirely formal. Whereas the set-theoretic semantic and the course-of-values semantic are extensional, the model-theoretic semantic is intensional, This is the contemporary successor of the distinction between an *Inhaltslogik* and the *Folgerungscalcul* about which Voigt, Husserl and their colleagues argued.

Associated with the *logica utens*/*logica docens* distinction is the logic as calculus/logic as language distinction. *Logic as calculus* is understood as a combinatorial tool for the formal manipulation of elements of a universe of discourse. Typically, but not necessarily, this universe of discourse is well-defined. We should, perhaps, better the logic as calculus on a purely syntactic level. The "Booleans" (and, although van Heijenoort did not specifically mention them, the Postulate theorists), reserved a formal deductive system for combinatorial-computational manipulation of the syntactic elements of their system. The semantic interpretation of the syntactic elements to be manipulated was external to the formal system itself. The semantic interpretation was given by the chosen universe of discourse. Again, the axioms selected for such a system were typically chosen to suit the needs of the particular field of mathematics being investigated, as were the primitives that provided the substance of the elements of the universe of discourse, whether sets, classes, or propositions, or points, lines, and planes, or tables, chairs and beer mugs. On the other hand, the *logica docens* is intended as an all-purpose formal logical system which is applicable regardless of the universe of discourse which provides the contents for its manipulation, regardless of the primitive terms upon which it operates, or their semantic interpretation, if any. It is in these terms that van Heijenoort also therefore distinguishes between *relativism* and *absolutism* in logic; a *logica docens* is appropriate relative to its specific universe of discourse; a *logica utens* is absolute in being appropriate to any and every universe of discourse. More broadly, there are many *logica utens*es, but only one, universally applicable, *logica docens*.

Pietarinen, for one, would agree with van Heijenoort with regard at least to Peirce, that his work belongs to logic as calculus; as Pietarinen expresses it [Pietarinen *2009b*, 19], "in relation to the familiar division between language as a universal medium of expression and language as a reinterpretable calculus," Peirce and his signifist followers took language to serve the latter role." Elsewhere, Pietarinen makes the case even more strongly and explicitly, asserting [Pietarinen *2009a*, 45] that: "Peirce's disaffection with unreasonably strong realist assumptions is shown by the fact than he did not advocate any definite, universal logic that would deserve the epithet of being the logic of our elementary thought. Logical systems are many, with variable interpretations to be used for the various purposes of scientific inquiry."

For Russell, as for Frege, says van Heijenoort, it is the character of this inclusiveness that makes their logical systems suitable not merely as a *calculus ratiocinator*, but as a *lingua characteristica* or *chrcteristica universalis*.[31] Thus, Frege's *Begriffsschrift* and Whitehead and Russell's *Principia Mathematica* are *both* calculus and language *at once*. Moreover, Frege would argue is unlike the calculi of the Booleans, not simply *both* calculus and language, but a language first and foremost. As we know, Schröder and Peano would argue over whether the classical Boole-Schröder or the logic of the *Formulaire* was the better pasigraphy, or *lingua universalis*, and Frege and Schröder, along the same lines, whether the *Begriffsschrift* or the classical Boole-Schröder was a *lingua*, properly so-called, and, if so, which was the better. Van Heijenoort would argue for the correctness of Frege's appraisal. In "Über die Begriffsschrift des Herrn Peano und meine einige" [Frege *1896*], Frege meanwhile argued that Peano's logical system, in the *Arithmetics principia* [Peano *1889*] and *Notations de la logique mathématique* [Peano *1894*], merely tended towards being a *characteristica* while yet remaining a *calculus ratiocinator*. Thus, Frege insists that only his Begriffsschrift is truly both a calculus and a language. He writes [*1896*, 371] (in van Heijenoort's translation [van Heijenoort *1967b*, 325, *n*. 3]): "Boole's logic is a *calculus ratiocinator*, but no *lingua charaterica*; Peano's mathematical logic is in the main a *lingua characterica* and subsidiarily, also a *calculus ratiocinator*, while my *Begriffsschrift* intends to be both with equal stress." Van Heijenoort [*1967b*, 325, *n*. 3] understands this to mean that Boole has a sentential calculus but no quantification theory; Peano has a notation for quantification theory but only a very deficient technique of derivation; Frege has a notation for quantification theory and a technique of derivation."

The discussion of the question of the nature and fruitfulness of a logical notation was an integral aspect of the debate as to whether a logical system was a calculus, a language, or both. In discussions on these issues with Frege and Peano, Schröder was the defender of the Peirce-Mitchell-Schröder pasigraphy against both Frege's *Begriffsschrift* and Peano's notation; see e.g. [Frege *1880/81; 1882*] for Frege's critique of Boole's algebra of logic as a mere calculus; [Frege *1895*; *1896*] for Frege's critiques of Schröder's algebra of logic and of Peano's axiomatic system respectively, [Schröder *1898a*; *1898b*] for a comparison of Peano's and Peirce's pasigraphies and defense of Peirce's pasigraphy, and [Peano *1898*] for Peano's reply to Schröder, for the major documents in this aspect of the discussion, and [Peckhaus *1990/1991*] for an account of Schröder's discussions with Peano. Peirce [*1906*] for his part regarded the logic of Peano's *Formulaire*, as presented in his [*1894*] *Notations de logique mathématique (Introduction au Formulaire de mathématiques)*, "no calculus; it is nothing but a pasigraphy…."

Returning to van Heijenoort's list of properties, it should be clear from the evidence which we have presented that, under this interpretation, Peirce indeed had both a *calculus ratiocinator*, or sentential calculus with derivation, defined in terms of illation (property 1) and a *characterica universalis*, or quantification theory and notation for quantification theory (property 3), and that these are clearly present in a single unified package in "On the Algebra of Logic: A Contribution to the Philosophy of Notation" [Peirce *1885*].

The other aspect of this universality is that, as a language, it is not restricted to a specific universe of discourse, but that it operates on the universal domain, what Frege called the *Universum*. Thus, the universe of discourse for Frege and Russell is the universal domain, or the universe. It is in virtue of the *Begriffsschrift*'s and the *Principia* system's universe of discourse being *the universe*, that enables these logical systems to say (to put it in colloquial terms) everything about everything in the universe. One might go even further, and with van Heijenoort understand that, ultimately, Frege was able to claim that there are only two objects in the *Universum*: *the True* and *the False*, and that every proposition in his system assigns the *Bedeutung* of a proposition to one or the other.

Johannes Lenhard [*2005*] reformulates van Heijenoort's distinction between logic as calculus and logic as language in ontological term, by suggesting that the concept of logic as a language upholds a model carrying an ontological commitment, and arguing in turn that Hilbert's formalism, expressed in his indifference to whether our axioms apply to points, lines, and planes or to tables, chairs, and beer mugs, bespeaks a model theory which is free of any ontology. It is precisely in this sense that Lenhard [*2005*, 99] cites Hintikka [*1997*, 109] as conceiving of Hilbert as opening the path to model theory in the twentieth century. This view of an ontologically challenged conception was already articulated by Hilbert's student Paul Isaac Bernays (1888–1977), who defined mathematical "existence" in terms of constructibility and non-contradictoriness within an axiom system and in his [*1950*] "Mathematische Existenz und Widerspruchsfreiheit".

What makes the logic of the *Begriffsschrift* (and of the *Principia*) a language preeminently, as well as a calculus, rather than a "mere" calculus, was that it is a *logica docens*, and it is absolute. The absoluteness guarantees that the language of the *Begriffsschrift* is a language, and in fact a universal language, and fulfills the Leibniz programme of establishing it as a *mathesis universa*, which is both a language and a calculus. In answer to the question of what Frege means when he says that his logical system, the *Begriffsschrift*, is like the language Leibniz sketched, a *lingua characteristica*, and not merely a logical calculus, [Korte *2010*, 183] says that: "According to the nineteenth century studies, Leibniz's *lingua characteristica* was supposed to be a language with which the truths of science and the constitution of its concepts could be accurately expressed." [Korte *2010*, 183] argues that "this is exactly what the *Begriffsschrift* is: it is a language, since, unlike calculi, its sentential expressions express truths, and it is a characteristic language, since the meaning of its complex expressions depend only on the meanings of their constituents and on the way they are put together." Korte argues that, contrary to Frege's claims, and those by van Heijenoort and Sluga in agreement with Frege, the *Begriffsschrift* is, indeed, a language, but *not* a calculus.[32]

Because of this universality, there is, van Heijenoort argues, nothing "outside" of the *Universum*. (This should perhaps set us in mind of Wittgenstein, and in particular of his proposition 5.5571 of the *Tractatus logico-philosophicus* [Wittgenstein *1922*], that "The limits of my language are the limits of my world"—"Die Grenzen meiner Sprache bedeuten die Grenzen meiner Welt.") If van Heijenoort had cared to do so, he would presumably have quoted Prop. 7 from the *Tractatus*, that, by virtue of the universality of the *logica docens* and its universal universe of discourse, anything that can be said must be said within and in terms of the *logica docens* (whether Frege's variant or Whitehead-Russell's), and any attempt to say anything *about* the system is "wovon man nicht sprechen kann." In van Heijenoort's terminology, then, given the universality of the universal universe of discourse, one cannot get outside of the system, and the system/metasystem distinction becomes meaningless, because there is, consequently, nothing outside of the system. It is in this respect, then, that van Heijenoort argued that Frege and Russell were unable to pose, let alone answer, metalogical questions about their logic. Or, as Wittgenstein stated it in his *Philosophische Grammatik* [Wittgenstein *1973*, 296]: "*Es gibt keine Metamathematik*," explaining that "Der Kalkül kann uns nicht prinzipielle Aufschlüssen über die Mathematik geben," and adding that "Es kann darum auch keine "führenden Probleme" der mathematischen Logik geben, denn das wären solche…."

It was, as van Heijenoort [*1967b*; *1977*; *1986b*; *1987*], Goldfarb [*1979*], and Gregory H. Moore [*1987*; *1988*] established, the model-theoretic turn, enabled by the work of Löwenheim, Skolem, and Herbrand, in turn based upon the classical Boole-Peirce-Schröder calculus, and opened the way to asking and treating metasystematic questions about the logical systems of Frege and Russell, and, as Moore [*1987*; *1988*] also showed, helped establish the first-order functional calculus of Frege and Russell as the epitome of "mathematical" logic.

Turning then specifically to Peirce, we can readily associate his concept of a *logica docens* as a general theory of semiotics with van Heijenoort's conception of Frege's *Begriffsschrift* and Whitehead-Russell's *Principia* as instances of a *logica docens* with logic as language; and likewise, we can associate Peirce's concept of *logica utens* as with van Heijenoort's concept of algebraic logic and the logic of relatives as instances of a *logica utens with* logic as a calculus. It is on this basis that van Heijenoort argued that, for the "Booleans" or algebraic logicians, Peirce included, the algebraic logic of the Booleans was *merely* a calculus, and not a language. By the same token the duality between the notions of logic as a calculus and logic as a language is tantamount to Peirce's narrow conception of logic as critic on the one hand and to his broad conception of logic as a general theory of signs or semiotics. It is on this basis that Volker Peckhaus has concluded [Peckhaus *1990/1991*, 174–175] that in fact Peirce's algebra and logic of relatives "wurde zum pasigraphischen Schlüssel zur Schaffung einer schon in den frühen zeichentheoretischen Schriften programmatisch geforderten wissenschaftlichen Universalsprache und zu einem Instrument für den Aufbau der "absoluten Algebra", einer allgemeinen Theorie der Verknüpfung," that is, served as both a *characteristica universalis* and as a *calculus ratiocinator*, the former serving as the theoretical foundation for the latter. Quoting Peirce from a manuscript of 1906 in which Peirce offered a summary of his thinking on logic as calculus and logic as language, Hawkins leads us to conclude that Peirce would not be content to consider satisfactory a logic which was merely a calculus, but not also a language, or pasigraphy also; Peirce, comparing his dual conception of logic as both calculus and language with the conceptions which he understood to be those of Peano on the one hand and of Russell on the other, writes in the manuscript "On the System of Existential Graphs Considered as an Instrument for the Investigation of Logic" of *circa* 1906 (notebook, MS. 499, 1-5, as quoted by [Hawkins *1997*, 120]):

> The majority of those writers who place a high value upon symbolic logic treat it as if its value consisted in its mathematical power as a calculus. In my [1901] article on the subject in Baldwin's Dictionary I have given my reasons for thinking…if it had to be so appraised, it could not be rated as much higher than puerile. Peano's system is no calculus; it is nothing but a pasigraphy; and while it is undoubtedly useful, if the user of it exercises a discreet [*sic*] freedom in introducing additional signs, few systems have been so wildly overrated as I intend to show when the second volume of Russell and Whitehead's Principles of Mathematics appears.[33] …As to the three modifications of Boole's algebra which are much in use, I invented these myself,—though I was anticipated [by De Morgan] as regards to one of them,—and my dated memoranda show…my aim was…to make the algebras as analytic of reasonings as possible and thus to make them capable of every kind of deductive reasoning…. It ought, therefore, to have been obvious in advance that an algebra such as I am aiming to construct could not have any particular merit [in reducing the number of processes, and in specializing the symbols] as a calculus.

Taking Peirce's words here at face value, we are led to conclude that, unlike those Booleans who were satisfied to devise calculi which were not also languages, Peirce, towards the conclusion of his life, if not much earlier, required the development of a logic which was both a calculus (or critic, "which are much in use, I invented these myself") and a language (or semiotic), and indeed in which the semiotic aspect was predominant and foundational, while considering the idea of logic as language as of paramount utility and importance.

We should also take into account that, for Frege, Dedekind and Russell, a salient feature of logicism as a philosophy was to understand logic to be the *foundation* of mathematics. That is, Frege's *Begriffsschrift* and *Gundlagen der Mathematik* and Whitehead and Russell's *Principia Mathematica* were designed firstly and foremostly, to carry out the construction of mathematics on the basis of a limited number of logical concepts and inference rules. Peirce agreed, in 1895 in "The Nature of Mathematics" (see [Peirce *2010*, 3–7]), that mathematics is a deductive science, arguing that it is also purely hypothetical, that [Peirce *2010*, 4] "a hypothesis, in so far as it is mathematical, is mere matter for deductive reasoning," and as such are imaginary—or ficitive—entities and thus not concerned with truth, and, moreover, not only stands on its own accord and as such is independent of logic (see, also [De Waal *2005*]), but argues also that mathematical reasoning does not require logic, while logic is exclusively the study of signs. This would appear to be as far as one could possibly get from the logicist position of constructing mathematics from the ground up from a small number of logical primitives and on the basis of one or a few rules of logical inference that in turn depend upon one or a small number of primitive logical connectives. Beginning with his father's assertion that mathematics is "the science which draws necessary conclusions," Peirce [*2010*, 7] supposed that his father most likely therefore understood it to be the business of the logician to formulate hypotheses on the basis of which conclusions are to be deduced. In response, he asserted that [Peirce *2010*, 7]:

> It cannot be denied that the two tasks, of framing hypotheses for deduction and of drawing the deductive conclusions are of widely different characters; nor that the former is similar to much of the work of the logician.

But there the similarity ends. Logic is, in his view, the science of semiotics. He explains that [Peirce *2010*, 7]:

> Logic is the science which examines signs, ascertains what is essential to being sign and describes their fundamentally different varieties, inquires into the general conditions of their truth, and states these with formal accuracy, and investigates the law of development of thought, accurately states it and enumerates its fundamentally different modes of working.

That being the case, it is difficult to see that Peirce would have adhered to a philosophy of logicism, in any case of the variety of logicism espoused by Frege or Russell. Cut it would very emphatically open the way to the use of logic for philosophical investigations (especially for philosophy of language), if it did not, indeed, equate logic, *quâ* semiotics, with philosophy of language and render it more than suitable, in that guise, for philosophical investigations. With this in mind, one might be reminded more of Rudolf Carnap's application of the "logische Analyse der Sprache", the logical analysis of language, for dissolving metaphysical confusions [Carnap *1931-32*], and of Wittgenstein's program of logico-linguistic therapy for philosophy, epitomized perhaps best by forays into logico-linguistic analysis such as those of J. L. Austin, Gilbert Ryle, and Peter Strawson, than of Frege's and Russell's efforts to construct all of mathematics on the basis of a second-order functional calculus or set theory, or even Carnap's *logische Aufbau der Welt* [Carnap *1928*] and Russell's [*1918*] "Philosophy of Logical Atomism". In case of doubt that Peirce would reject the primacy, or foundational aspect, of logic over mathematics, he continued in "The Nature of Mathematics" [Peirce *2010*, 8], to assert that, although it "may be objected" that his just-stated definition

> places mathematics above logic as a more abstract science the different steps [of] which must precede those of logic; while on the contrary logic is requisite for the business of drawing necessary conclusions.

His immediate response [Peirce *2010*, 8] is: "But I deny this." And he explains [Peirce *2010*, 8]:

> An application [of] logical theory is only required by way of exception in reasoning, and not at all in mathematical deduction.

This view is indubitably the explanation for what auditors at a meeting of the New York Mathematical Society (forerunner of the American Mathematical Society) in the early 1890s must have considered Peirce's emotional outburst, recorded for posterity by Thomas Scott Fiske (1865–1944) who adjudged Peirce as having a "dramatic manner" and "reckless disregard of accuracy" in cases dealing with what Peirce considered to be "unimportant details," describing ([Fiske *1939*, 15; *1988*, 16]; see also [Eisele *1988*]) "an eloquent outburst on the nature of mathematics" when Peirce

> proclaimed that the intellectual powers essential to the mathematician were "concentration, imagination, and generalization." Then, after a dramatic pause, he cried: "Did I hear some one say demonstration? Why, my friends," he continued, "demonstration is merely the pavement upon which the chariot of the mathematician rolls."

On this basis, we are, I might suggest, safe in supposing that, while Peirce accepted and indeed promoted logic as a significant linguistic tool, in the guise of semiotics, for the investigation and clarification of philosophical concepts, we are equally safe in supposing that he did not likewise accept or promote the foundational philosophy of logicism in any manner resembling the versions formulated by Frege, Dedekind, or Russell. Moreover, this interpretation is affirmed in Peirce's explicit account of the relation between mathematics, logic, and philosophy, especially metaphysics, found in *The Monist* article "The Regenerated Logic" of 1896; there, Peirce [*1896*, 22–23] avows the "double assertion, first that logic ought to draw upon mathematics for control of disputed principles," rather than building mathematics on the basis of principles of logic, "and second that ontological philosophy ought in like manner to draw upon logic," since, in the hierarchy of abstractness, mathematics is superior to logic and logic is superior to ontology. Peirce, in "The Simplest Mathematics", explicitly rejected Dedekind's assertion, from *Was sind und was sollen die Zahlen?* [Dedekind *1888*, VIII] that arithmetic (algebra and analysis) are branches of logic (see [Peirce *2010*, 32–33]), reiterating the hypothetical nature of mathematics and the categorical nature of logic, and reinforcing this point by adding the points that, on the contrary, formal logic, or mathematical logic, is mathematics, indeed depends upon mathematics, and that formal logic is, after all, only a small part of logic, and, moreover, associating of logic with ethics, as a normative science, dependent upon ethics more even than upon mathematics.

Finally, returning to a theme emphasized by van Heijenoort in close conjunction with his distinction between logic as calculus and logic as language, there is the distinction between the semantic and the syntactic approaches to logic.

As [Zeman *1977*, 30] remarks, Peirce's

> basic orientation toward deductive logic is a *semantical* one, as we might be led to expect from his association of "logic proper" with the *object* of a sign. The icons of the algebra of logic are justified by him on what we recognize as truth-functional, and so semantic, grounds (see [Peirce *1933a*, 3.38 *f.*], for example) and the most basic sign of the systems of existential graphs, the "Sheet of assertion" on which logical transformations are carried out, is "considered as representing the universe discourse" [Peirce *1933b*, 4.396]; such representation is a semantical matter. But contemporary logic makes a distinction that Peirce did not make. It is necessary to stud logic not only from a radically semantical point of view, in which propositions are thought of as being true or false, but also from a *syntactic* or *proof theoretical* point of view, in which the deducibility of propositions from each other is studied without reference to interpretations in universes of any sort and so without reference to truth and falsity.
>
> Peirce failed to distinguish between logic as proof-theoretical and logic as semantical, but he can hardly be faulted for that; Gottlob Frege, who with Peirce must be considered a co-founder of contemporary logic, also failed to make the distinction,[34] and even Whitehead and Russell are fuzzy about it. Indeed, a clear recognition of the necessity for distinguishing between logical syntax and semantics does not arise until later, with the developments in logic and the foundations of math which culminated in Gödel's celebrated completeness and incompleteness results of 1930 and 1931 respectively.[35]

For us, and for van Heijenoort, the syntax of a logic consists of the uninterpreted symbols of the system, together with its inference rules. The semantics is comprised of the universe or universes of discourse which the represent, that is, the (extra-logical) interpretation that we give to the symbols of the system, whether, in the words of Hilbert, these are points, lines, and planes or tables, chairs, and beer mugs. For Russell, the semantic is a set-theoretic one, composed of Cantor's sets; for Frege, it is a course-of-values (*Werthverläufe*) semantic. In these terms, the principal difference between the algebraic logicians, including Peirce, on the one hand, and the "logisticians", Frege and Russell in particular, as well as Hilbert, is that the universe or universal domain, Frege's *Universum*, serves for the semantic interpretation of the system, whereas, for the algebraic logicians, following De Morgan (who coined the concept of the "universe of discourse" as under the term "universe of a proposition, or of a name" [De Morgan *1846*, 380; *1966*, 2]), it is a restricted and pre-determined subset or the *Universum* or universal domain, or, in van Heijenoort's modern terms, a *model*. A model is simply an ordered pair comprised of the elements of a (non-empty) domain ($D$) and one of its subdomains ($S$), and a mapping from the one to the other, in which the mapping assigns symbols to the elements of the domain and subdomain. When we change the semantic interpretation for the domain, say from points, lines, and planes, to tables, chairs, and beer mugs, the obtain a new model, logically equivalent to to the old model, and truth is preserved from the one model to the other. If every such model has the same cardinality, then is logically equivalent, and we cannot distinguish one from the other interpretation, and if truth is preserved for all of the interpretations for every formula in these models, we say that the models are *categorical*. In terms of our consideration, it was van Heijenoort's claim that the algebraic logicians, Peirce among them, devised logical calculi, but not languages, because, unlike Frege and Russell, theirs was a model-theoretic, or intensional, rather than a set-theoretic, or extensional, semantic.

7. *Peirce's distinguishing singular propositions, such as "Socrates is mortal" from universal propositions such as "All Greeks are mortal":*

Some difficulties with the conditions for validity in traditional syllogistic logic that arose when compared with conditions for validity in George Boole's treatment of the syllogism in his algebraic logic led directly to Peirce's contributions to the development of algebraic logic. In his first publication in logic, Peirce dealt with a means of handling singular propositions, which, in traditional logic, had most typically been treated as special cases of universal propositions. Peirce's efforts to distinguish singular from universal propositions led him to his first work in developing a theory of quantification for Boole's algebraic logic. The problem of distinguishing singular from universal propositions was it is fair to suggest, a major motivation for Peirce in devising his algebraic logic as an improvement over that of Boole, if not, indeed, the sole motivation. It is clear, we should also be reminded, that with the introduction of empty classes in Boole's system, some syllogisms that were valid in traditional, *i.e*. Aristotelian, logic, which assumes that there are no nonempty classes, were not necessarily valid in the algebraic logic, while some syllogisms that are not valid in the algebra of logic are valid in traditional logic.[36] The issue of existential import arises when dealing with singular propositions because traditional logic assumes the existence or non-emptiness of referents of classes for universal propositions but, in treating singular propositions as special cases of universal propositions, likewise assumes the existence of referents singular terms.[37]

We must readily grant that the issue of the distinction between singular and universal propositions was raised by Frege and taken up in earnest by Peano and Russell. We already remarked, in particular, on Frege's [*1893*] *Grundgesetze der Arithmetik*, especially [Frege *1893*, §11], where Frege's function \ξ replaces the definite article, such that, for example, \ξ(positive √2) represents the concept which is the proper name of the positive square root of 2 when the value of the function \ξ is the positive square root of 2, and on Peano's [*1897*] "Studii di logica matematica", in which Peano first considered the role of "the", the possibility its elimination from his logical system; whether it can be eliminated from mathematical logic, and if so, how. In the course of these discussions, Russell raised this issue with Norbert Wiener (see [Grattan-Guinness *1975*, 110]), explaining that:

> There is need of a notation for "the". What is alleged does not enable you to put "0 = etc. Df". It was a discussion on this very point between Schröder and Peano in 1900 at Paris that first led me to think Peano superior.

But the significance of this distinction between singular and universal propositions, which he learned from the medieval logicians, who, however, were unable to hand it satisfactory, was an early impetus for Peirce when he began his work in logic.[38]

Peirce recognized that Boole judged the adequacy of his algebra of logic as a system of general logic by how it compared with syllogistic, and he himself nevertheless sought to provide a syllogistic pedigree for his algebra of logic, he rejected the notion that the algebra of logic is a conservative extension of syllogistic. Thus, in the four page "Doctrine of Conversion" (1860) from loose sheets of Peirce's logic notebook of 1860–1867 [Peirce *1860-1867*], Peirce was rewriting the syllogism "All men are animals. *X* is a man. Therefore *X* is an animal" as "If man, then animal. But man. Therefore animal."

Some of Peirce's earliest work in algebraic logic, and his earliest published work, therefore was undertaken precisely for the purpose of correcting this defect in Boole's logical system, in particular in "On an Improvement in Boole's Calculus of Logic" [Peirce *1868*]. This effort can be traced to at least 1867. The goal of improving Boole's algebraic logic by developing a quantification theory which would introduce a more perspicacious and efficacious use of the universal and existential quantifiers into Boole's algebra and likewise permit a clear distinction between singular propositions and universal propositions as early as 1867, it was not until 1885 that Peirce was able to present a fully-articulated first-order calculus along with a tentative second-order calculus. He was able to do so on the basis of notational improvements that were developed in 1883 by his student Oscar Howard Mitchell, whereby universal quan-tifiers are represented as logical products and existential quantifiers by logical sums in which both relations or "logical polynomials" and quantifiers are indexed. The indexing of quantifiers was due to Peirce; the indexing of the logical polynomials was due to Mitchell. Peirce then denoted the existential and universal quantifiers by '$\Sigma_i$' and '$\Pi_i$' respectively, as logical sums and products, and individual variables, *i*, *j*, …, are assigned both to quantifiers and predicates. He then wrote '$l_{i,j}$' for '*i* is the lover of *j*'. Then "Everybody loves somebody" is written in Peirce's quantified logic of relatives as $\Pi_i \Sigma_j l_{i,j}$, i.e. as "Everybody is the lover of somebody". In Peirce's own exact expression, as found in his "On the Logic of Relatives" [*1883*, 200], we have: "$\Pi_i \Sigma_j l_{i,j} > 0$ means that everything is a lover of something." The quantifiers run over classes, whose elements are itemized in the sums and products, and an individual for Peirce is one of the members of these classes.

Despite this, singular propositions, especially those in which definite descriptions rather than proper names, have also been termed "Russellian propositions", so called because of their designation by Bertrand Russell in terms of the iota quantifier or iota operator, employing an inverted iota to be read "the individual *x*"; thus, e.g. $(\iota x)\Phi(x)$.[39] In *Principia Mathematica* [Whitehead & Russell *1910*, 54], Russell writes "$\phi!x$" for the first-order function of an individual, that is, for any value for any variable which involves only individuals; thus, for example, we might write $\mu!$(Socrates) for "Socrates is a man". In the section on "Descriptions" of *Principia* [Whitehead & Russell *1910*, 180], the iota operator replaces the notation "$\phi!x$" for singulars with "$(\iota x)(\Phi x)$" so that one can deal with definite descriptions in as well as names of individuals. This is a great simplification of the lengthy and convoluted explanation that Russell undertook in his *Principles of Mathematics* [Russell *1903*, 77–81, §§76–79], in order to summarize the points that: '"Socrates is a-man" expresses identity between Socrates and one of the terms denoted by *a man*' and '"Socrates is one among men," a proposition which raises difficulties owing to the plurality of men' in order to distinguish a particular individual which is an element of a class having more than one member from a unary class, *i.e.* a class having only one member.

§3. *How much, if anything Peirce and Frege knew of, or influenced, the other's work*. In the "Preface" to the English translation of the *Algebra of Logic* of Louis Couturat (1868–1914) Philip Edward Bertrand Jourdain (1879–1919) summarized what he considered to be the true relation between the algebraic logicians on the one hand and the logisticians, Frege, Peano, and Russell on the other, writing [Jourdain *1914*, VIII]:

> We can shortly, but very fairly accurately, characterize the dual development of the theory of symbolic logic during the last sixty years as follows: The *calculus ratiocinator* aspect of symbolic logic was developed by BOOLE, DE MORGAN, JEVONS, VENN, C. S. PEIRCE, SCHRÖDER, Mrs. LADD-FRANKLIN and others; the *lingua characteristica* aspect was developed by FREGE, PEANO and RUSSELL. Of course there is no hard and fast boundary-line between the domains of these two parties. Thus PEIRCE and SCHRÖDER early began to work at the foundations of arithmetic with the help of the calculus of relations; and thus they did not consider the logical calculus merely as an interesting branch of algebra. Then PEANO paid particular attention to the calculative aspect of his symbolism. FREGE has remarked that his own symbolism is meant to be a *calculus ratiocinator* as well as a *lingua characteristica*, but the using of FREGE's symbolism as a calculus would be rather like using a three-legged stand-camera for what is called "snap-shot" photography, and one of the outwardly most noticeable things about RUSSELL's work is his combination of the symbolisms of FREGE and PEANO in such a way as to preserve nearly all of the merits of each.

As we in turn survey the distinguishing characteristics of mathematical logic and, equivalently, enumerate the innovations which Frege was the first to devise and, in so doing, to create mathematical logic, we are led to conclude that many, if not all, of these were also devised, to greater or lesser extent, also Peirce, if not entirely simultaneously, then within close to the same chronological framework. One major difference is that Frege is credited, for the greater part correctly, with publishing nearly all of the elements within one major work, the *Begriffsschrift* of 1879, whereas Peirce worked out these elements piecemeal, over a period that began in or around 1868 and over the next two decades and beyond, in widely scattered publications and unpublished manuscripts. From Peirce's biography, we can extract the explanation; that his time and efforts were distracted in many directions, owing in no small measure to his personal circumstances, and to a significant degree as a result of his lack of a long-term academic position.

Beyond this explanation, we are led to also inquire to what extent, if any Peirce was influenced by Frege, and also whether Frege was influenced by Peirce, and, in both cases, if so, to what extent. So far as I am aware, the most recent and thoroughgoing efforts to deal with these questions are to be found in [Hawkins *1971*; *1993*; *1997*, 134–137]. It is therefore useful to apply the evidence that [Hawkins *1971*; *1993*; *1997*, 134–137] provides in order to understand how "Peircean" was the "Fregean" revolution.

The short answer is that there is strong circumstantial evidence that Peirce and his students were well aware of the existence of Frege's *Begriffsschrift* beginning at least from Peirce's receipt of Schröder's [*1880*] review of that work, but it is uncertain whether, and if so how deeply, Peirce or his students studied the *Begriffsschrift* itself; but most likely they viewed it from Schröder's [*1880*] perspective. It is unclear whether Schröder and met Peirce met during Peirce's sojourn in Germany in 1883 (see [Hawkins *1993*, 378–379, *n*. 3]); neither is it known whether Schröder directly handed the review to Peirce, or it was among other offprints that Schröder mailed to Peirce in the course of their correspondence (see [Hawkins *1993*, 378]). We know at most with certainty that Ladd-Franklin [*1883*, 70–71] listed it as a reference in her contribution to Peirce's [*1883a*] *Studies in Logic*, but did not write about it there, and that both Allan Marquand (1853–1924) and the Johns Hopkins University library each owned a copy; see [Hawkins *1993*, 380]), as did Christine Ladd-Franklin (see [Anellis *2004-2005*, 69, *n*. 2]. The copy at Johns Hopkins has the acquisition date of 5 April 1881, while Peirce was on the faculty. It is also extremely doubtful that Peirce, reading Russell's [*1903*] *Principle of Mathematics* in preparation for a review of that work (see [Peirce *1966*, VIII, 131, *n*. 1]), would have missed the many references to Frege. On the other hand, the absence of references to Frege by Peirce suggests to Hawkins [*1997*, 136] that Peirce was unfamiliar with Frege's work.[40] What is more likely is that Peirce and his students, although they knew of the existence of Frege's work, gave it scant, if any attention, not unlike the many contemporary logicians who, attending to the reviews of Frege's works, largely dismissive when not entirely negative, did not study Frege's publications with any great depth or attention.[41]

In virtue of Peirce's deep interest in and commitment to graphical representations of logical relations, as abundantly evidenced by his entitative and existential graphs, and more generally, in questions regarding notation and the use of signs, he might well have been expected to develop an interest in and study of Frege's three-dimensional notation (see, *e.g*. **Figure 3** for Frege's *Begriffsschrift* notation and proof that $(\forall x)(f(x) = g(x))$).

On the other side of the equation, there is virtually no evidence at all that Frege had any direct, and very little indirect, knowledge of the work of Peirce beyond what he found in the references to Peirce in Schröder's *Vorlesungen über die Algebra der Logik*. As [Hawkins *1993*, 378] expresses it: "It is almost inconceivable that Schröder, in his [*1880*] review of Frege's *Begriffsschrift*, with his references to four of Peirce's papers, would have escaped Frege's notice, with Schroder so subject as he is, to Frege's critical ire" (see [Dipert *1990/1991a*, 124–125]). But this hardly entails any serious influence of Peirce upon Frege's work beyond a most likely dismissive attitude that would derive from the derivative rejection of Schröder's work on Frege's part, devolving as a corollary upon the work as well of Peirce. [Dipert *1990/1991*, 135] goes further and is far less equivocal, remarking that there was "no discernible influence of Frege" on the work even of German logicians or set theorists, and barely any of Russell.[42] If, therefore, we take literally our question of how Peircean was the Fregean revolution, our reply must be: "Not at all," and certainly in no obvious wise. As [Hawkins *1993*, 379] concludes, "Frege's writings appear to be quite unaffected by Peirce's work (see [Hawkins *1975*, 112, *n*. 2; *1981*, 387])." Neither has any evidence been located indicating any direct communication between Peirce and Frege, or between any of Peirce's students and Frege, while relations between Frege and Schröder were critical on Schröder's side towards Frege, and hostile on Frege' side towards Schröder, without, however, mentioning him by name (see [Dipert *1990/1991*, 124–125]). What can be said with surety is that Frege left no published evidence of references to Peirce's work. Whether there ever existed unpublished pointers to Peirce, either in Frege's writings or in his library, we can only guess inasmuch as his *Nachlaß* has not been seen since the bombing during World War II of the University of Münster library where it had been deposited for safekeeping.[43]

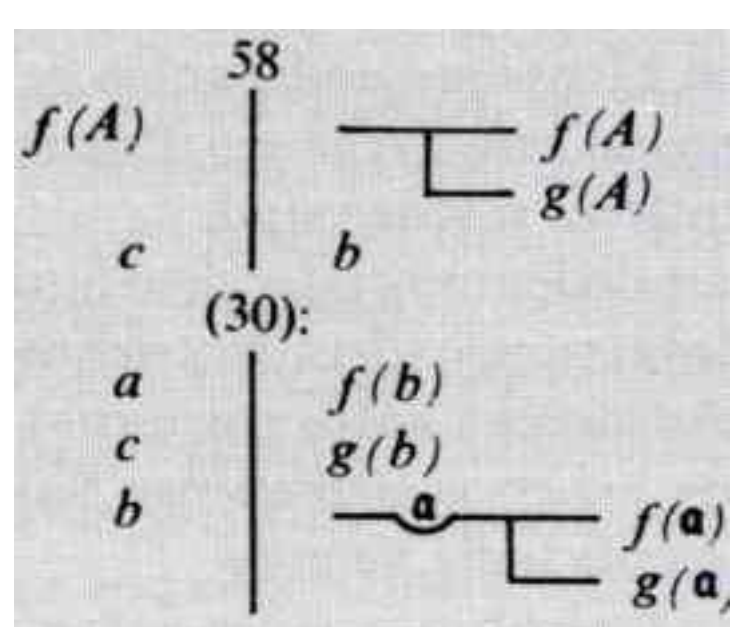

**Figure 3. Frege's *Begriffsschrift* notation and proof that $(\forall x)(f(x) = g(x))$**

The comparison is especially helpful when we consider the same logical formula as expressed in the notations, set side-by-side of Frege, Peirce's graphical notation, and Peano-Russell (in the modified version presented in Whitehead and Russell's *Principia Mathematica*, and in the Polish notation. For this purpose, we turn to [Leclercq *2011*], which presents a table (**Table 4**) comparing the representation of the formula [(~*c* ⊃ *a*) ⊃ (~*a* ⊃ *c*)] ⊃ {(~*c* ⊃ *a*) ⊃ [(*c* ⊃ *a*) ⊃ *a*]} in Frege's notation, Peirce's graphical notation, the Peano-Russell notation, and the Polish notation of Łukasiewicz, as follows:

| | |
|---|---|
| Algebraic: | Peano-Russell: [(~*c* ⊃ *a*) ⊃ (~*a* ⊃ *c*)] ⊃ {(~*c* ⊃ *a*) ⊃ [(*c* ⊃ *a*) ⊃ *a*]}<br>Peirce: [($\bar{c}$ —< *a*) —< ( $\bar{a}$ —< *c*)] —< {($\bar{c}$ —< *a*) —< [(*c* —< *a*) —< *a*]}<br>Schröder: [(*c*´ € *a*) € (*a*´ € *c*)] € {(*c*´ € *a*) € [(*c* € *a*) € *a*]}<br>Łukasiewicz: CCCN*ca*CN*ac*CCN*ca*CC*caa* |
| "graphical": | **Frege :** 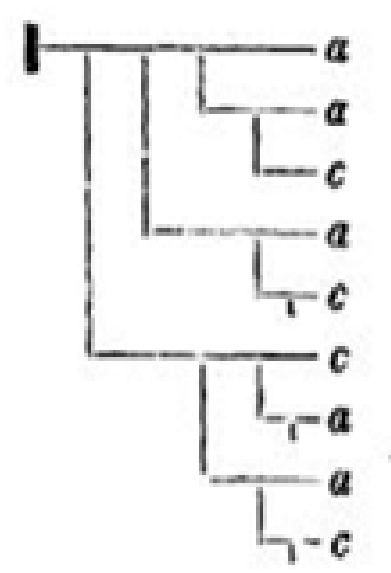 **Peirce :** 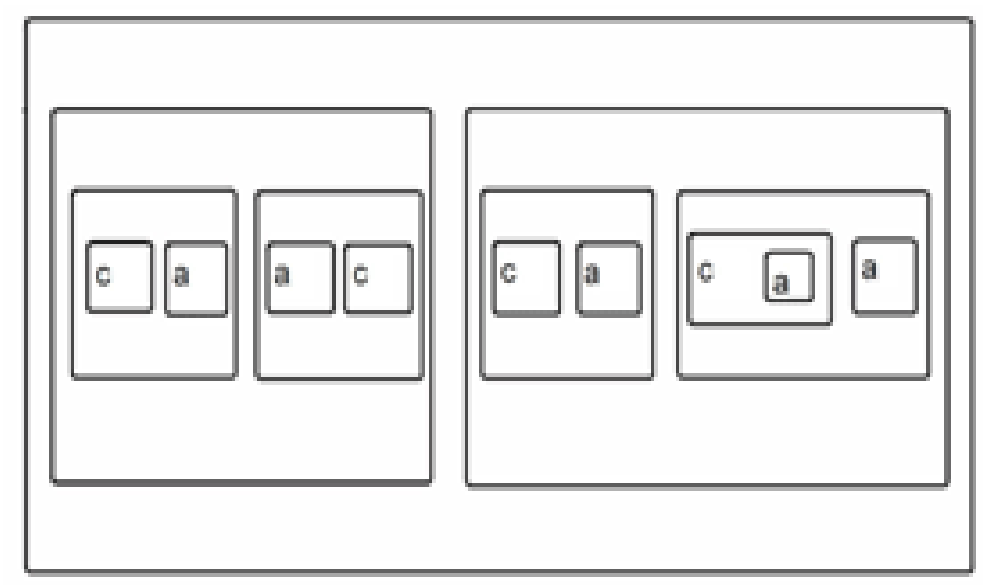 |

**Table 4. Algebraic and Graphical Notations for [(~*c* ⊃ *a*) ⊃ (~*a* ⊃ *c*)] ⊃ {(~*c* ⊃ *a*) ⊃ [(*c* ⊃ *a*) ⊃ *a*]}**

On the other side of the equation, there is virtually no evidence at all that Frege had any direct, and very little indirect, knowledge of the work of Peirce beyond what he found in the references to Peirce in Schröder's *Vorlesungen über die Algebra der Logik*. As [Hawkins *1993*, 378] expresses it: "It is almost inconceivable that Schröder, in his [*1880*] review of Frege's *Begriffsschrift,* with his references to four of Peirce's papers, would have escaped Frege's notice, with Schroder so subject as he is, to Frege's critical ire" (see [Dipert *1990/1991a*, 124–125]). But this hardly entails any serious influence of Peirce upon Frege's work beyond a most likely dismissive attitude that would derive from the derivative rejection of Schröder's work on Frege's part, devolving as a corollary upon the work as well of Peirce. [Dipert *1990/1991*, 135] goes further and is far less equivocal, remarking that there was "no discernible influence of Frege" on the work even of German logicians or set theorists, and barely any of Russell.[44] If, therefore, we take literally our question of how Peircean was the Fregean revolution, our reply must be: "Not at all," and certainly in no obvious wise. As [Hawkins *1993*, 379] concludes, "Frege's writings appear to be quite unaffected by Peirce's work (see [Hawkins *1975*, 112, *n*. 2; *1981*, 387])." Neither has any evidence been located indicating any direct communication between Peirce and Frege, or between any of Peirce's students and Frege, while relations between Frege and Schröder were critical on Schröder's side towards Frege, and hostile on Frege' side towards Schröder, without, however, mentioning him by name (see [Dipert *1990/1991*, 124–125]). What can be said with surety is that Frege left no published evidence of references to Peirce's work. Whether there ever existed unpublished pointers to Peirce, either in Frege's writings or in his library, we can only guess inasmuch as his *Nachlaß* has not been seen since the bombing during World War II of the University of Münster library where it had been deposited for safekeeping.[45]

So far as concerns Russell's knowledge of the work in logic of Peirce and Schröder, Danielle Macbeth [*2005*, 6] is willing to admit only that "Russell did not learn quantificational logic from Frege. But he did learn a logic of relations from Peirce," and, referring specifically to Peirce's [*1870*] and [*1883*], moreover that "he knew Peirce's work while he was developing the polyadic predicate calculus, and in particular, he knew that Peirce had developed a complete logic of relations already in 1883 based on Boole's logical algebra…."[46] In fact, Russell read, an composed notes on, Peirce's [*1880*] "On the Algebra of Logic" and [*1885*] "On the Algebra of Logic: A Contribution to the Philosophy of Notation" [Russell *ca. 1900-1901*] while undertaking his own earliest compositions on the logic of relations and claiming to have himself invented the logic of relations (see [Anellis *1995*]). He at the same time read and composed notes [Russell *1901b*] and marginal annotations on Schröder's *Vorelsungen über die Algebra der Logik*, which the inscriptions indicate he acquired in September 1900; he also owned an offprints of Schröder's [*1877*] *Der Operationskreis des Logikkalkül* and [*1901*] "Sur une extension de l'idée d'ordre" (see [Anellis *1900/1901*]).

Being as precise and cautious as possible, we can assert confidently that Frege did not undertake to examine Schröder's work until faced with Schröder's review of the *Begriffsschrift*, did not undertake to examine Boole's work until chided by Schröder and Venn for ignoring it in the *Begriffsschrift*, and, if there exists any documentary evidence that Frege read any of Peirce's work, it has not yet come to light, whereas we have excellent documentation of what Russell read and when he read it, not only in his notes on Peirce and Schröder, but also, for some of his readings, a log recording titles and approximate dates.[47] All of the published reviews of Frege's *Begriffsschrift* comment upon the difficulty and "cumbrousness" of his notation and the total absence of reference or mention, direct or indirect, to any contemporary work in logic. It was precisely the second complaint, in particular the one emanating from Schröder, that impelled Frege to undertake to compare his work with that of Boole, beginning with the unpublished "Booles rechnende Logik und die Begriffsschrift" [Frege *1880/81*] and "Booles logische Formelsprache und die Begriffsschrift" [Frege *1882*].

In both directions whether Frege to Peirce or Peirce to Frege, there is, therefore, no direct reference to the writings of the one in the writings of the other.

I am also well aware that there are some other historians of logic who seek to undertake to redress the balance between the "Booleans" and the "Fregeans", among them Randall Dipert, but who also caution against zealotry in the opposite direction which credits Peirce too greatly as having anticipated Frege. In the case, for example, of propositional logic and the development of truth tables, Dipert [*1981*] notes that Peirce (1) was not averse to employing numerical values, not just 1 and 0, for evaluating the truth of propositions, rather than *true* and *false*, but a range of numerical values, and (2) that the formulas employed employed by Peirce were, depending upon the context in which they occurred, allowed to have different interpretations, so that their terms might represent classes rather than propositions; and hence it would be over-simplifying the history of logic to argue that Peirce was a precursor in these respects of Frege, or anticipated Frege or some one else. In this Dipert is essentially reiterating one of Husserl's arguments in support of his assertion [Husserl *1891b*, 267] that Schröder's algebra of logic is a calculus but not a language, and hence not a logic properly so-called, assigning values of 1 and 0 rather than true and false, to his equations. Taking Dipert's caution into consideration, this should not, however, detain the historian from acknowledging Peirce for his accomplishments.

If, therefore, we interpret our original question to mean: to what extent did Peirce (and his students and adherents) obtain those elements that characterize the "Fregean" revolution in logic?, our reply must be: "To a considerable extent," but not necessarily all at once and in one particular publication. Nor were they themselves always fully cognizant of distinctions which we would today insist upon and retrospectively either impute to them or blame them for not recognizing. It also helps if we keep in mind the simple truth that, whether one is speaking of Peirce or Frege, or for that matter any historical development, simplifications much more often than not obscure the true picture of that development. As the Maschellin says in "Der Rosenkavalier", "In dem 'Wie', da liegt der ganze Unterschied." In particular, we cannot pinpoint an exact date or even specific year for the "birth" of "modern" or "mathematical" logic, remembering that Peirce, Frege, and their colleagues were working out the details of their logical innovations over the course of time, piecemeal.

In the words of Jay Zeman: with respect to the technical aspects of logic [Zeman *1986*, 1], if we want to provide a nice summary without arguing over priority or picking dates: "Peirce developed independently of the Frege-Peano-Russell (FPR) tradition all of the key formal logical results of that tradition. He did this first in an algebraic format similar to that employed later in *Principia Mathematica*…."

*Notes*

[1] Husserl sent Frege offprints of both his [Husserl *1891a*] and [Husserl *1891b*]; see [Gabriel 1*980*, 30]; and as [Pietersma *1967*, 298, *n*. 1] also notes, it was the review of Schröder that induced Frege to write his [*1895*] critique of Schröder.

[2] See [Hamacher-Hermes *1991*] for details.

[3] [Voigt *1892*, 292]: "Die Algebra der Logik will jedoch mehr als dieses Zugeständniss; sie will ein Logik in vollen Sinne sein, sie behauptet wesentlich denselben Inhalt zu haben, dieselben Ziele zu verfolgen wie di ältere Logik und das zwar auf einem sichereren, exacteren Weg."

[4] Russell's notes on [Peirce *1880*] and [Peirce *1885*] (ms., 3pp.) date from *ca*. 1900-1901; his notes on Schröder, *Vorlesungen über die Algebra der Logik*, ms. 6pp., RA file #230: 030460, date from 1901. They are lodged in the Bertrand Russell Archive, Ready Memorial Library, McMaster University, Hamilton, Ontario, Canada.

[5] See [Russell *1903*, 10*n*, 12*n*, 13, 22, 24, 26, 142, 201*n*, 221*n*, 232, 306*n*, 320*n*, 367*n*] for mentions of Schröder and [Russell *1903*, 23, 26, 203*n*, 232*n*, 320*n*, 376, 387*n*] for even fewer mentions of Peirce.

[6] [Goldfarb *1979*], [Moore *1987*; *1988*], and [Peckhaus *1992*; *1994*] examine those technical elements of the history of logic in the 1920s that helped defined the canonical conception of mathematical logic as, first and foremost, first-order predicate logic. [Anellis *2011*] surveys the conflicting views on logic during this period by those who were working in logic at that time, and in particular contemporary attitudes towards Russell's place in the then-current development of logic.

[7] Peano made this point himself quite clearly in a letter to Russell of 19 March 1901 (quoted in [Kennedy *1975*, 206]), declaring that Russell's paper on the logic of relations [Russell *1901*] "fills a gap between the work of Peirce and Schröder on the one hand and the *Formulaire* on the other."

In the "Preface of "On Cardinal Numbers", Whitehead [*1902*, 367] explains that the first two sections serve as "explaining the elements of Peano's developments of mathematical logic…and of Russell's symbolism for the Logic of Relations" as given in [Russell *1901*] and expresses his belief that "these two methods are almost indispensible for the development of the theory of Cardinal Numbers." Section III, he notes [Whitehead *1902*, 368], on "Finite and Infinite Cardinal Numbers" [Whitehead *1902*, 378–383], was "entirely due to Russell and is written by him throughout." Regarding the Peano notation as modified by Russell, Whitehead [*1902*, 367] judges it its invention "forms an epoch in mathematical reasoning."

[8] Held in Box 2 of Van Heijenoort Nachlaß: Papers, 1946–1983; Archives of American Mathematics, University Archives, Barker Texas History Center, University of Texas at Austin.

[9] Peirce's Nachlass originally located in Harvard University's Widener Library and since located in Harvard's Houghton Library, with copies of all materials located in the Max H. Fisch Library at the Institute for American Thought, Indiana University-Purdue University at Indianapolis [IUPUI].

[10] [Bernatskii *1986; 1990*] considers Peirce on the Aristotelian syllogism.

[11] [Shosky *1997*, 12, *n*. 6] cites Post's [*1921*] in [van Heijenoort *1967a*, 264–283]; but see also the doctoral thesis [Post *1920*]. Łukasiewicz is mentioned, but [Shosky *1997*] gives no reference; see [Łukasiewicz *1920*]. See also [Anellis *2004a*].

[12] Zhegalkin employed a technique resembling those employed by Peirce-Mitchell-Schröder, Löwenheim, Skolem, and Herbrand to write out an expansion of logical polynomials and assigning them Boolean values.

[13] [Shosky *1997*] totally ignores the detailed and complex history of the truth-table method and shows no knowledge of the existence of the truth-table device of Peirce in 1902. For a critique of Shosky's [*1997*] account and the historical background to Peirce's work, see [Anellis *2004a*; *forthcoming (b)*]; see also [Clark *1997*] and [Zellweger *1997*] for their re-(dis)-covery and exposition of Peirce's work on truth-functional analysis and the development of his truth-functional matrices.

[Grattan-Guinness *2004-05*, 187–188], meanwhile, totally misrepresents the account in [Anellis *2004a*] of the contributions of Peirce and his cohorts to the evolution and development of truth-functional analysis and truth tables in suggesting that: (1) Peirce did not achieve truth table matrices and (2) that [Anellis *2004a*] was in some manner attempting to suggest that Russell somehow got the idea of truth tables from Peirce. The latter is actually strongly contraindicated on the basis of the evidence that was provided in [Anellis *2004a*], where it is shown that Peirce's matrices appeared in unpublished manuscripts which did not arrive at Harvard until the start of 1915, after Russell had departed, and were likely at that point unknown, so that, *even if* Russell could had been made aware of them, it would have more than likely have been from Maurice Henry Sheffer (1882–1964), and *after* 1914.

[14] Peirce's orginal tables from MS 399 are reproduced as plates 1-3 at [Fisch & Turquette *1966*, 73–75].

[15] In *Logic (Logic Notebook 1865–1909)*; MS 339:440–441, 443; see [Peirce *1849–1914*], which had been examined by Fisch and Turquette. On the basis of this work, Fisch and Turquette [*1966*, p. 72] concluded that by 23 February 1909 Peirce was able to extend his truth-theoretic matrices to three-valued logic, there-by anticipating both Jan Łukasiewicz in "O logice trójwartosciowej" [Łukasiewicz *1921*], and Emil Leon Post in "Introduction to a General Theory of Elementary Propositions" (Post *1921*), by a decade in developing the truth table device for triadic logic and multiple-valued logics respectively. Peirce's tables from MS 399 are reproduced as the plates at [Fisch & Turquette *1966*, 73–75].

[16] See, *e.g*. [Pycior *1981*; *1983*] on the algebraic predisposition of Boole, De Morgan, and their colleagues and on De Morgan's work in algebra respectively; [Laita *1975*; *1977*] on the role played for Boole by analysis in the inception and conception of logic, and [Rosser *1955*] on Boole on functions.

[17] *Elective symbols*, *x*, *y*, *z*, etc., are so called in virtue of the nature of the operation which they are understood to represent, expressions which involve these are called *elective functions*, equations which contain elective functions are called *elective equations*; and an elective operation or function on *xy* is one which picks out, or selects", in succession, those elements of the class *Y* which are also members of the class *X*, hence successively selects all (and only) members of both *X* and *Y* (see [Boole *1847*, 16]).

[18] The equivalence was in principle apparently recognized both by Frege and Russell, although they worked in opposite directions, Frege by making functions paramount and reducing relations as functions, Russell by silently instituting the algebraic logician's rendition of functions as relations and, like De Morgan, Peirce and Schröder, making relations paramount. Thus, in [Oppenheimer & Zalta *2011*, 351], we read: "Though Frege was interested primarily in reducing mathematics to logic, he succeeded in

reducing an important part of logic to mathematics by defining relations in terms of functions. In contrast, Whitehead and Russell reduced an important part of mathematics to logic by defining functions in terms of relations (using the definite description operator).We argue that there is a reason to prefer Whitehead and Russell's reduction of functions to relations over Frege's reduction of relations to functions. There is an interesting system having a logic that can be properly characterized in relational type theory (RTT) but not in functional type theory (FTT). This shows that RTT is more general than FTT. The simplification offered by Church in his FTT is an over-simplification: one cannot assimilate predication to functional application."

[19] Russell [*1903*, 187] evidently bases his complaint upon his between "Universal Mathematics", meaning universal algebra in what he understood to be Whitehead's sense, and the "Logical Calculus", the former "more formal" than the latter. As Russell understood the difference here, the signs of operations in universal mathematics are variables, whereas, for the Logical Calculus, as for every other branch of mathematics, the signs of operations have a constant meaning. Russell's criticism carries little weight here, however, inasmuch as he interprets Peirce's notation as conflating class inclusion with set membership, not as concerning class inclusion and implication. In fact, Peirce intended to deliberately allow his "claw" of illation () to be a primitive relation, subject to various interpretations, including, among others, the ordering relation, and material implication, as well as class inclusion, and—much later— set membership.

[20] See [Merrill *1978*] for a study of De Morgan's influences on Peirce and a comparison of De Morgan's and Peirce's logic of relations.

[21] This list was not meant to be exhaustive, let alone current.

[22] On Boole's influence on Peirce, see, *e.g.* [Michael *1979*].

[23] Aspects of these developments were also considered, e.g. in [van Heijenoort *1982*; *1986*, 99–121], [Anellis *1991*], [Moore *1992*].

[24] As [Anellis *1992*, 89–90] first suggested, there is "a suggestion of relativization of quantifiers is detectable in the *Studies in logical algebra* (MS 519) dating from May of 1885 and contemporaneous with Peirce's vol. 7 AJM paper of 1885 *On the algebra of logic: a contribution to the philosophy of notation,* where we can find something very much resembling what today we call *Skolem normal form*; [and] in *The logic of relatives: qualitative and quantitative* (MS 532) dating from 1886 and scheduled for publication in volume 5, we may look forward to Peirce's use of a finite version of the Löwenheim-Skolem theorem, accompanied by a proof which is similar to Löwenheim's" (see [Peirce *1993*, 464, *n*. 374.31–36]).

[25] Among the studies of Peirce's work on set theory and the logic of numbers in additions to [Shields *1981*; *1997*] and [Gana *1985*] are [Dauben *1977*; *1981*], [Levy *1986*], [Myrvold *1995*], and [Lewis *2000*].

[26] Among those far better equipped than I to entertain questions about the history of the philosophical aspects of this particular discussion and in particular whether or not Peirce can be regarded as a logicist; see, *e.g.* [Houser *1993*], part of a lengthy and continuing debate; see also, e.g. [Haack *1993*], [Nubiola *1996*], and [De Waal *2005*]. This list is far from exhaustive.

[27] In a much broader sense, Peirce would also include induction and deduction as means of inference along with deduction.

[28] In addition to [van Heijenoort *1967b*], [van Heijenoort *1977*; *1986b*; *1987*] are also very relevant to elucidating the distinctions logic as calculus/logic as language, *logica utens*/*logica docens*, and relativism/absolutism in logic and their relation to the differences between "Booleans" and "Fregeans". The theme of logic as calculus and logic as language is also taken up, *e.g.*, in [Hintikka *1988*; *1997*].

[29] In its original conception, as explicated by the medieval philosophers, the *logica utens* was a practical logic for reasoning in specific cases, and the *logica docens* a teaching logic, or theory of logic, concerning the general principles of reasoning. These characterizations have been traced back at least to *the Logica Albertici Perutilis Logica* of Albertus de Saxonia (1316–1390) and his school in the 15th century, although the actual distinction can be traced back to the *Summulae de dialectica* of Johannes Buridanus (*ca*. 1295 or 1300–1358 or 1360). See, *e.g.* [Bíard *1989*] for Buridan's distinction, and [Ebbesen *1991*] and [Bíard *1991*] on Albert. The distinction was then borrowed by Peirce; see [Bergman & Paavola *2003-*] on "Logica docens" and "Logica utens".

[30] The concept of universe of discourse originated with Augustus De Morgan in his "On the Syllogism, I: On the Structure of the Syllogism" [*1846*, 380; *1966*, 2], terming it the "universe of a proposition, or of a name" that, unlike the fixed universe of all things that was employed by Aristotle and the medieval logicians, and remained typical of the traditional logic, "may be limited in any manner expressed or understood." Schröder came subsequently to conceive this in terms of a *Denkbereich*. When considering, then, Schröder's terminology, one is dealing with a different *Gebietkalkül* for the different *Denkbereichen*.

[31] See, *e.g*. [Patzig *1969*] for an account of Frege's and Leibniz's respective conceptions of the *lingua characteristica* (or *lingua charactera*) and their relationship. [Patzig *1969*, 103] notes that Frege wrote of the idea of a *lingua characteristica* along with *calculus ratiocinator*, using the term "lingua characteristica" for the first time only in 1882 in "Über den Zweck der Begriffsschrift" (in print in [Frege *1883*]) and then again in "Über die Begriffsschrift des Herrn Peano und meine einige" [Frege *1897*], in the *Begiffsschrift* [Frege *1879*] terming it a "formal language"—"Formelsprache". In the foreword to the *Begriffsschrift*, Frege wrote only of a *lingua characteristica* or "allgemeine Charakteristik" and a *calculus ratiocinator*, but not of a *charcteristica universalis*.

[32] Sluga [*1987*, 95] disagrees with van Heijenoort only in arguing the relevance of the question of the claimed universality for the logic of the *Begriffsschrift*; in Frege, Sluga argues, concepts result from analyzing judgments, whereas in Boole, judgments results from analyzing concepts. [Korte *2010*, 285] holds, contrary to all other suggestions on the question of why the *Begriffsschrift* properly is a language, that it is so because of Frege's logicism, and nothing else.

In attempting to understand Peirce's position with regard to whether his logic of relatives is a language or a calculus, we are therefore returned to the debate regarding the question of whether, and if so, how and to what extent, Peirce was a logicist in the sense of Frege and Peano; see *supra*, *n*. 18.

[33] Peirce tended to conflate Russell and Whitehead even with respect to Russell's *Principles of Mathematics*, even prior to the appearance of the co-authored *Principia Mathematica* [Whitehead and Russell *1910-13*], presumably because of their earlier joint work "On Cardinal Numbers" [Whitehead *1902*] in the *American Journal of Mathematics*, to which Russell contributed the section on, a work with which Peirce was already familiar.

[34] [Zerman *1977*, 30] refers, imprecisely, to the English translation at [van Heijenoort *1967a*, 13] of [Frege *1879*], but not to Frege's original, and does not explain the nature of Frege's failure to distinguish logic as proof-theoretic from semantic. The reference is to the discussion in [Frege *1879*, §4].

[35] [Zeman *1977*, 30] has in mind [Gödel *193*0; *1931*].

[36] The problem of the existential import thus became a major issue of contention between traditional logicians and mathematical or symbolic logicians; [Wu *1962*] is a brief history of the question of existential import from Boole through Strawson.

[37] I would suggest that it is in measure precisely against this backdrop that Russell became so exercised over Meinongian propositions referring to a non-existent present bald king of France as to not only work out his theory of definite descriptions in "On Denoting" [Russell *1905*], but to treat the problem in *The Principles of Mathematics* [Russell *1903*] and in *Principia Mathematica* [Whitehead & Russell *1910*], as a legitimate problem for logic as well as for philosophy of language and for metaphysics.

[38] See, *e.g*. [Michael *1976a*] on Peirce's study of the medieval scholastic logicians. [Martin *1979*] and [Michael *1979b*] treat individuals in Peirce's logic, and in particular in Peirce's quantification theory. [Dipert *1981*, 574] mentions singular propositions in passing in the discussion of Peirce's propositional logic.

[39] In a letter to Jourdain of March 15 1906, Russell wrote (as quoted in [Russell *1994*, xxxiii]: "In April 1904 I began working at the Contradiction again, and continued at it, with few intermissions, till January 1905. I was throughout much occupied by the question of Denoting, which I thought was probably relevant, as it proved to be. ...The first thing I discovered in 1904 was that the variable denoting function is to be deduced from the variable propositional function, and is not to be taken as an indefinable. I tried to do without *i* as an indefinable, but failed; my success later, in the article "On Denoting", was the source of all my subsequent progress."

[40] We cannot entirely discount the relevance of Peirce's isolation from his academic colleagues following his involuntary resignation from his position at Johns Hopkins University in 1884 (see, e.g. [Houser *1990/1991*, 207]). This may perhaps underline the significance of Schröder's seeking, in his letter to Peirce of March 2, 1897 (see [Houser *1990/1991*, 223], to "to draw your attention to the pasigraphic movement *in Italy.* Have you ever noticed the 5 vols. of Peano's Rivista di Matematica together with his "Formulario" and so many papers of Burali-Forti, Peano, Pieri, de Amicis, Vivanti, Vailati, etc. therein, as well as in the reports of the Academia di Torino and in other Italian periodicals (also in the German Mathematische Annalen)."

[41] Reviews of the *Begriffsschrift* have been collected and published in English in [Frege *1972*]; the introductory essay [Bynum *1972*] presents the familiar standard interpretation that Frege's work was largely ignored until interest in it was created by Bertrand Russell; more recent treatments, e.g.[Vilkko *1998*], argue that Frege's work was not consigned to total oblivion during Frege's lifetime.

[42] [Gray *2008*, 209] would appear to suggest a strong influence "of the algebraic logic of Boole, Jevons, and others" on Frege, who thought that "a deep immersion [in their work] might drive logicians to new and valuable questions about logic." There is no virtually evidence for accepting this supposition, however.

[43] See [Wehmeier & Schmidt am Busch *2000*] on the fate of Frege's *Nachlaβ*.

[44] [Gray *2008*, 209] would appear to suggest a strong influence "of the algebraic logic of Boole, Jevons, and others" on Frege, who thought that "a deep immersion [in their work] might drive logicians to new and valuable questions about logic." There is no virtually evidence for accepting this supposition, however.

[45] See [Wehmeier & Schmidt am Busch *2000*] on the fate of Frege's *Nachlaβ*.

[46] [Macbeth *2005*] is clearly unaware of [Russell *ca. 1900-1901*], Russell's notes on [Peirce *1880*] and [Peirce *1885*].

[47] See [Russell *1983*] for Russell's log "What Shall I Read?" covering the years 1891–1902.